\documentclass[11pt]{article}
\usepackage{amsmath}
\usepackage{amssymb}
\usepackage{amsbsy}
\usepackage{amsthm}
\usepackage{epsfig}
\usepackage{wrapfig}
\usepackage{eepic}
\usepackage{color}
\usepackage{graphicx}
\usepackage{amsfonts}%

\newcommand{\bu}{{\bf u}}

\newcommand{\be}{{\bf e}}
\newcommand{\bD}{{\bf D}}
\newcommand{\bU}{{\bf U}}

\newcommand{\bv}{{\bf v}}
\newcommand{\bw}{{\bf w}}
\newcommand{\bx}{{\bf x}}

\newcommand{\bH}{{\bf H}}
\newcommand{\bX}{{\bf X}}

\newcommand{\bphi}{{\boldsymbol \phi}}

\def\div{\operatorname{div}}

\def\curl{\operatorname{curl}}
\newcommand{\blf} {\mathbf{f}}
\newcommand{\dO} {{\partial\Omega}}

\begin{document}

\title{Efficient discretizations for the EMAC formulation of the incompressible Navier-Stokes equations}

\author{
Sergey Charnyi
\footnote{Department of Mathematical Sciences, Clemson University, Clemson, SC 29634 (scharny@clemson.edu).}
\and
Timo Heister
\footnote{Department of Mathematical Sciences, Clemson University, Clemson, SC 29634 (heister@clemson.edu), partially supported by NSF Grant DMS1522191.}
\and
Maxim A. Olshanskii
\footnote{Department of Mathematics, University of Houston, Houston TX 77004 (molshan@math.uh.edu), partially supported by Army Grant 65294-MA and NSF Grant DMS1522252.}
\and
Leo G. Rebholz\footnote{Department of Mathematical Sciences, Clemson University, Clemson, SC 29634 (rebholz@clemson.edu), partially supported by Army Grant 65294-MA and NSF Grant DMS1522191.}
}
\date{}

\maketitle

\begin{abstract}
We study discretizations of the incompressible Navier-Stokes equations, written in the newly developed energy-momentum-angular momentum conserving (EMAC) formulation.  We consider linearizations of the problem, which at each time step will reduce the computational cost, but can alter the conservation properties.  We show that a skew-symmetrized linearization delivers the correct balance of (only) energy and that the Newton linearization conserves momentum and angular momentum, but conserves energy only up to the nonlinear residual.  Numerical tests show that linearizing with 2 Newton steps at each time step is very effective at preserving all conservation laws at once, and giving accurate answers on long time intervals.  The tests also show that the skew-symmetrized linearization is significantly less accurate.
The tests also show that the  Newton linearization of EMAC finite element formulation compares favorably to other traditionally used finite element formulation of the   incompressible Navier-Stokes equations in primitive variables.
\end{abstract}

\section{Introduction}

We consider discretizations of the incompressible Navier-Stokes equations (NSE), which are given by
\begin{eqnarray}
\bu_t + (\bu\cdot\nabla) \bu + \nabla p - \nu\Delta \bu & = & \blf, \label{nse1} \\
\div\bu & = & 0,\\
\bu(0) &= &\bu_0,
\end{eqnarray}
in a domain $\Omega\subset \mathbb{R}^d$, $d$=2 or 3, {\color{black}and for $t>0$}, where $\bu$ and $p$ represent velocity and pressure, $\blf$ is the external forcing, $\bu_0$ represents the initial velocity, and $\nu$ the kinematic viscosity.  The system must also be equipped with appropriate boundary conditions, see e.g. \cite{temam}.  This system models the evolution of water, oil, and air flow (air under 220 m.p.h.), and therefore is important in a wide array of problems in science and engineering.

In the recent work \cite{CHOR17}, the authors showed that if the NSE system was discretized by a Galerkin method using the reformulation
\[
\bu\cdot\nabla \bu + \nabla p = 2\bD(\bu)\bu + (\div \bu)\bu + \nabla P,
\]
with $P=p-\frac12 |\bu|^2$ and $\bD$ denoting the rate of deformation tensor, then each of energy, momentum, angular-momentum, 2D enstrophy, helicity, and total vorticity would all be correctly balanced, {\it even if the method does not enforce the divergence constraint strongly.}  For this reason, we call the `$2\bD(\bu)\bu + (\div \bu)\bu$' formulation of the nonlinearity the EMAC formulation (energy, momentum, angular momentum conserving).  In most common Galerkin methods for incompressible flow problems, such as mixed finite element methods, the divergence constraint is only enforced weakly \cite{JLMNR17}, and in \cite{CHOR17} we show how this is the main cause of conservation law violation if standard formulations of the nonlinearity are used.  However, with the EMAC formulation, all of these conservation laws are obeyed in Galerkin discretizations. The overall efficiency of the formulation and it clear superior performance over traditionally used  finite element formulations for the NSE in velocity--pressure variables (for flow examples where conservation properties are of importance) have been recently reported in  \cite{CHOR17,SL17}.

The purpose of this paper is to investigate linearizations of the EMAC formulation at each time step in a temporal discretization.  It is common practice in simulations of the NSE to approximate the discrete system at each time step, by approximating the nonlinear term with some linear approximation, e.g. for convective form, replacing $\bu^n \cdot\nabla \bu^n$ at time step $n$ by $\bU \cdot\nabla \bu^n$, where $\bU$ is some known approximation of $\bu^n$ such as $2\bu^{n-1}-\bu^{n-2}$.  The advantage of such a linearization is that now the system requires just one linear solve per time step, instead of needing multiple linear solves to resolve a Newton or Picard iteration, but it retains energy stability.  However, such modifications of the nonlinearity can adversely affect conservation properties of a numerical scheme.  We will consider below linearizations of the EMAC nonlinearity.  First, we derive a skew-symmetric form of the EMAC nonlinearity, and show it conserves energy but not momentum or angular momentum.  We then derive the Newton linearization, and show it preserves momentum and angular momentum conservation in the EMAC nonlinearity.  We also study the energy balance of the Newton linearization and conclude the energy is balanced correctly up to certain terms, which vanish if one let the Newton method converge. Otherwise (if, for example,
only one iteration of the Newton method is done), these terms are of higher order compared to approximation error of the method.
After analyzing these linearizations, we provide results for numerical tests that show the skew-symmetric linearization provides rather inaccurate solutions, but the Newton linearization (with 2 or less iterations per time step) gives very good results.

This paper is arranged as follows.  In section 2, we give mathematical preliminaries, and define notation, to allow for a smooth presentation in later sections.  Section 3 derives and analyzes the two linearizations of the EMAC nonlinearity which are discussed above.  Section 4 numerically tests these linearizations on benchmark problems, and finally conclusions and future directions are discussed in section 5.

\section{Notation and preliminaries}

We consider the domain $\Omega \subset \mathbb{R}^d$, d=2 or 3,  and denote the $L^2(\Omega)$ inner product and norm on $\Omega$ by $(\cdot,\cdot)$ and $\| \cdot \|$, respectively.  The natural velocity and pressure spaces for the NSE are
\begin{eqnarray*}
\bH^1_0(\Omega) = \{ \bv \in H^1(\Omega)^d,\ \bv |_{\partial\Omega}=0 \}, \qquad
L^2_0(\Omega) = \{ q \in L^2(\Omega),\ \int_{\Omega} q\ dx=0 \}.
\end{eqnarray*}

We will consider subspaces $\bX \subset \bH^1_0(\Omega)$, $Q \subset L^2_0(\Omega)$ to be finite dimensional, and more specifically that $\bX$ and $Q$ are finite element velocity and pressure spaces corresponding to an admissible triangulation of $\Omega$.  We further assume, for simplicity, that $\bX$ and $Q$ satisfy inf-sup compatibility conditions~\cite{GR86}; (non inf-sup stable pairs require stabilization terms that will affect conservation properties, and should be studied case-by-case).  Our analysis can be easily extended to other types of Galerkin methods.

Most common finite element discretizations of Navier-Stokes and related systems, e.g. using Taylor-Hood elements, only enforce the divergence-free constraint $\div \bu=0$ weakly.  Instead of the pointwise constraint, a numerical solution $\bu $ from $\bX $ is enforced to satisfy
\[
(\div \bu ,q )=0 \quad \forall~ q \in Q.
\]
Even though convergence theory of mixed finite element methods guarantees $\| \div \bu \|$ be small in some sense (converges to 0 with optimal spatial rate), it is shown analytically and computationally in \cite{CHOR17} that it is large enough to cause physically-conserved quantities such as energy, momentum and angular momentum to not be conserved by the discretization of common NSE formulations like convective, skew-symmetric, and rotation forms.  Enlarging the pressure space $Q $ to enforce $\div\bX \subset Q $, which would provide pointwise enforcement of the divergence constraint, is usually not possible as it would violate the inf-sup compatibility condition and make the method numerically unstable (except in a few exceptional cases).


\subsection{Vector identities}
Consider $\bu,\ \bv,\ \bw\in \bH^1(\Omega)$,  and note that we do not enforce that any of these quantities are solenoidal.
Define the trilinear form $b:\bH^1(\Omega) \times \bH^1(\Omega) \times \bH^1(\Omega) \rightarrow \mathbb{R}$ by
\begin{equation}\label{conv_b}
b(\bu,\bv,\bw) = (\bu\cdot\nabla \bv,\bw).
\end{equation}
We recall the following properties of $b$.  The first two follow immediately from integration by parts, provided that the normal component of $\bu$ vanishes on $\dO$:
\begin{align}
b(\bu,\bv,\bw)& = -b(\bu,\bw,\bv)-( (\div \bu)\bv,\bw),  \label{vecid1} \\
b(\bu,\bw,\bw)& = -\frac12 \left( (\div \bu)\bw,\bw \right), \label{vecid2}\\
b(\bu,\bv,\bw)& = ( (\nabla \bv) \bu,\bw)=( (\nabla \bv)^T \bw,\bu).  \label{vecid3}
\end{align}
We denote the symmetric part of $\nabla \bu$ by
$
\nabla_s \bu :=  \bD(\bu) = \frac{ \nabla \bu + (\nabla \bu)^T }{2},
$
and the skew-symmetric part by
$
\nabla_n \bu := \frac{ \nabla \bu - (\nabla \bu)^T }{2}.
$
For any $\bu, \bv \in \bH^1(\Omega)$ one readily checks
\begin{equation}
(\nabla_n \bu)\bv =  \frac12 ( \curl{} \bu) \times \bv. \label{vecid6}
\end{equation}
Note that we define $\curl{} \bu$ in 2d in the usual way, as the 3d curl of $\bu$ extended by 0 in the
third component.

Straight-forward calculations and \eqref{vecid6} provide the following  vector identities for functions $\bu, \bv \in \bH^1(\Omega)$:
\begin{align}
(\bu\cdot \nabla) \bu & =  (\curl \bu) \times  \bu + \nabla \frac12 | \bu |^2 =: (\curl \bu) \times  \bu + \nabla q, \label{vecid5} \\
(\bu\cdot \nabla) \bu & =(\nabla \bu)\bu  =  (\nabla_s \bu)\bu + (\nabla_n \bu)\bu = \bD(\bu)\bu +  \frac12 (\curl \bu) \times \bu, \label{vecid6b}
\end{align}
where $q:= \frac{ | \bu |^2}{2}$.
Also note that identity \eqref{vecid6b} implies that
\begin{equation}
( \bD(\bu) \bu,\bu) = (( \nabla \bu)\bu,\bu) = b(\bu,\bu,\bu).\label{vecid6c}
\end{equation}
From \eqref{vecid6}--\eqref{vecid6b} we obtain the following representation of the inertia term from the momentum equations:
\begin{equation}
(\bu\cdot \nabla) \bu = 2\bD(\bu)\bu - \nabla q. \label{vecid7}
\end{equation}
The identity \eqref{vecid7} is key to the EMAC  formulation.

Using \eqref{vecid6c} and then \eqref{vecid2}, we obtain
\begin{equation}
 2(\bD(\bu )\bu ,\bu ) +  ((\div \bu )\bu ,\bu ) =   2b(\bu ,\bu ,\bu ) +  ((\div \bu )\bu ,\bu )=0. \label{econs}
\end{equation}

\subsection{Energy balance and conservation of linear and angular  momentum by NSE}

Assuming smoothness of solutions, the NSE solutions can be shown to deliver energy balance and conserve  linear momentum, and angular momentum, which we define by
\begin{eqnarray*}
&\text{Kinetic energy}\quad & E=\frac12(\bu,\bu):= \frac12 \int_\Omega |\bu|^2\mbox{d}\bx;\\
&\text{Linear momentum}\quad& M:=\int_\Omega \bu\, \mbox{d}\bx;\\
&\text{Angular momentum}\quad& M_\bx:=\int_\Omega \bu\times\bx\, \mbox{d}\bx.
\end{eqnarray*}
To see the balances, assume for simplicity that the solution $\bu,\ p$ have compact support in $\Omega$ (e.g. consider an isolated vortex), and test the NSE with $\bu$, $\be_i$ (the $i^{th}$ standard basis vector), and $\bphi_i =\be_i \times \bx$ to obtain
\begin{eqnarray*}
\frac{d}{dt}E + \nu \| \nabla \bu \|^2 & = & (\blf,\bu), \\
\frac{d}{dt} M_i& = & \int_{\Omega} f_i \ d\bx, \\
\frac{d}{dt} (M_\bx)_i & = & \int_{\Omega} (\blf\times \bx)_i\ d\bx,
\end{eqnarray*}
noting that each nonlinear and pressure term vanished, and using
\[
(\bu_t,\bphi_i)=\frac{d}{dt} \int_{\Omega} (\bu \times \bx)\cdot \be_i \ d\bx = \frac{d}{dt} (M_\bx)_i.
\]
For a numerical scheme to have physical accuracy, its solutions should admit balances that match these balances as close as
possible.  The key point here is that the nonlinear term does not contribute to any of these balances.

\section{Analysis of linearizations of the EMAC nonlinearity}

We will consider the following discretization of the NSE, for which the temporal discretization is left general:  At each time step, find $(\bu^n,p^n)\in (\bX,Q)$ satisfying
\begin{align}
(d_t (\bu^n),\bv) + 2(\bD(\bu^n)\bu^n,\bv) + ((\div \bu^n)\bu^n,\bv) - (p^n,\div \bv) + \nu (\nabla \bu^n,\nabla \bv) & =  (\blf^n,\bv), \label{emac1} \\
(\div \bu^n,q) & = 0, \label{emac2}
\end{align}
for all $(\bv,q)\in (\bX,Q)$.  Defining the time derivative term appropriately can produce, e.g., a BDF-k temporal discretization, and also a Crank-Nicolson temporal discretization if the $n$ superscript instead denotes a half step in time.

We consider now the  energy balance and the conservation of momentum, and angular momentum for the EMAC scheme, after carefully defining the problem setting.  Most physical boundary conditions will alter the physical balances of these quantities, as they should, when walls and interfaces are present.  Moreover, numerical treatment of boundaries (e.g. enforcing conditions strongly or weakly) will also affect the balances.
In this work, we  isolate the affect of the nonlinearity treatment on the quantities of interest from the contribution of the boundary conditions by assuming that the finite element solution $\bu^n$ and $p^n$ at each time step, along with the forcing $f^n$, is supported in a subset $\widehat{\Omega}\varsubsetneq \Omega$, i.e., there is a strip $S=\Omega \setminus \widehat{\Omega}$ along $\partial\Omega$ where each $\bu^n$ is zero.   Note that this implies there is a strip of elements along $\partial\Omega$  where $\bu^n$ and $p^n$  vanish.  The physical scenario for these boundary conditions  is  the evolution of an isolated vortex in a self-induced flow.  However, if a formulation/scheme is not able to conserve quantities $E,\ M,\ M_x$ in this setting, then the balances of these quantities in more general settings will also be incorrect and non-physical.

To derive the balances for energy, momentum, and angular momentum for the EMAC scheme above, we choose $\bv=\bu^n$, $\widetilde{\be_i}$, and $\widetilde{\bphi_i}$, respectively, where the tilde operator alters these quantities only in the strip of elements along the boundary so that they vanish on the boundary (note that $\widetilde{\be_i},\ \widetilde{\bphi_i} \in \bX$).  This gives

\begin{equation}
(d_t \bu^n,\bu^n) + (2\bD(\bu^n)\bu^n,\bu^n) + ((\div \bu^n)\bu^n,\bu^n ) + \nu \| \nabla \bu^n  \|^2 = (\blf^n,\bu^n ),
\end{equation}
\begin{equation}
\int_{\Omega} d_t u_i^n + (2\bD(\bu^n)\bu^n,\be_i) + ((\div \bu^n)\bu^n,\be_i )  = \int_{\Omega} f_i^n\ d\bx
\end{equation}
\begin{equation}
\int_{\Omega} d_t (\bu^n \times \bx)_i\ d\bx  + (2\bD(\bu^n)\bu^n,\bphi_i) + ((\div \bu^n)\bu^n,\bphi_i )  = \int_{\Omega} (\blf^n \times \bx)_i \ d\bx.
\end{equation}
Supposing that the temporal discretization is Crank-Nicolson, we obtain the time derivative terms
\begin{eqnarray*}
(d_t \bu^n,\bu^n) & = &  \frac{E(t^n) - E(t^{n-1})}{\Delta t}, \\
\int_{\Omega} d_t u_i^n & = & \frac{M_i(t^n) - M_i(t^{n-1})}{\Delta t}, \\
\int_{\Omega} d_t (\bu^n \times \bx)_i\ d\bx  & = & \frac{(M_x)_i(t^n) - (M_x)_i(t^{n-1})}{\Delta t}.
\end{eqnarray*}
Hence the key to balances of these quantities that are analogous to the true physical balances is that the nonlinear terms must vanish.  From \eqref{econs}, we obtain that
\[
(2\bD(\bu^n)\bu^n,\bu^n) + ((\div \bu^n)\bu^n,\bu^n ) =0.
\]
By expanding the rate of deformation tensor and using
$(\bu \cdot\nabla \bu ,\be_i) =  - ( (\div \bu )\bu ,\be_i)$ and then \eqref{vecid3}, we find that
\begin{eqnarray*}
 2(\bD(\bu^n )\bu^n ,\widetilde{\be_i}) + ((\div \bu^n )\bu^n ,\widetilde{\be_i}) &=& 2(\bD(\bu^n )\bu^n ,\be_i) + ((\div \bu^n )\bu^n ,\be_i)\\
 & = &  b(\bu^n,\bu^n ,\be_i) + b(\be_i,\bu^n,\bu^n) + ((\div \bu^n )\bu^n ,\be_i) \\
 & = &  b(\be_i,\bu^n,\bu^n)= (\be_i\times\bu^n,\bu^n)\\
 & = & 0
\end{eqnarray*}
since $\be_i$ is divergence-free.  Similarly,
\begin{eqnarray*}
2(\bD(\bu^n )\bu^n ,\widetilde{\bphi_i}) + ((\div \bu^n )\bu^n ,\widetilde{\bphi_i})&=&2(\bD(\bu^n )\bu^n ,\bphi_i) + ((\div \bu^n )\bu^n ,\bphi_i)\\
& = & b(\bu^n,\bu^n ,\bphi_i) + b(\bphi_i,\bu^n,\bu^n) + ((\div \bu^n )\bu^n ,\bphi_i)  \\
& = &   b(\bu^n ,\bu^n ,\bphi_i) + ((\div \bu^n )\bu^n ,\bphi_i) \\
& =& -b(\bu^n,\bphi_i,\bu^n)\\
& = & 0,
\end{eqnarray*}
with the last step following from expanding the expression and performing a simple calculation.  Thus we obtain the following balances for energy, momentum, and angular momentum,
\begin{equation}
 \frac{E(t^n) - E(t^{n-1})}{\Delta t} + \nu \| \nabla \bu^n  \|^2 = (\blf^n,\bu^n ),
\end{equation}
\begin{equation}
\frac{M_i(t^n) - M_i(t^{n-1})}{\Delta t}  = \int_{\Omega} f_i^n\ d\bx
\end{equation}
\begin{equation}
 \frac{(M_x)_i(t^n) - (M_x)_i(t^{n-1})}{\Delta t}  = \int_{\Omega} (\blf^n \times \bx)_i \ d\bx,
\end{equation}
which are discrete analogues for energy, momentum and angular momentum conservation.

We consider below these conservation properties for alterations of the EMAC scheme, where the nonlinear term is linearized.

\subsection{An energy stable linearization for EMAC}

By expanding the deformation tensor and using \eqref{vecid1}, we find that for $\bw,\bv\in \bX$,
\begin{eqnarray*}
2(\bD(\bw )\bw ,\bv) + ((\div \bw )\bw ,\bv)
& = &
b(\bw,\bw,\bv) + b(\bv,\bw,\bw) +  ((\div \bw )\bw ,\bv) \\
& = &
b(\bv,\bw,\bw) - b(\bw,\bv,\bw).
\end{eqnarray*}
This motivates the following EMAC linearization,
\begin{eqnarray}
(d_t (\bu^n),\bv) + b(\bv,\bu^n,\bu ^*) - b(\bu^n,\bv,\bu ^*)  - (p^n,\div \bv) + \nu (\nabla \bu^n,\nabla \bv) & = & (\blf^n,\bv), \label{emacL1}\\
(\div \bu^n,q) & = & 0, \label{emacL2}
\end{eqnarray}
where $\bu^*$ is a known approximation of $\bu(t)$.  For example, a first order approximation such as $\bu^*=\bu^{n-1}$ would be appropriate for backward Euler, $\bu^*=2\bu^{n-1} - \bu^{n-2}$ for BDF2, and $\bu^*=\frac32 \bu^{n-1} - \frac12 \bu^{n-2}$ for Crank-Nicolson (for this case the $n$ represents a half time step).

Since this linearized scheme matches \eqref{emac1}-\eqref{emac2} except for the nonlinear term, we need only consider the nonlinear term to study conservation properties of \eqref{emacL1}-\eqref{emacL2}.  Following the analysis above, we observe that conservation of energy, momentum and angular momentum will be attained only if the term $b(\bv,\bu^n,\bu ^*) - b(\bu^n,\bv,\bu ^*)$ vanishes for $\bv=\bu^n,\ \widetilde{\be_i},\ \widetilde{\bphi_i}$, respectively.

Energy conservation is achieved by construction, as we observe that for $\bv=\bu^n$,
\[
b(\bu^n,\bu^n,\bu ^*) - b(\bu^n,\bu^n,\bu ^*)=0.
\]
However, neither momentum or angular momentum are conserved by \eqref{emacL1}-\eqref{emacL2}, since
\begin{eqnarray*}
 b(\widetilde{\be_i},\bu^n,\bu ^*) - b(\bu^n,\widetilde{\be_i},\bu ^*)
 & = & b(\be_i,\bu^n,\bu ^*) - b(\bu^n,\be_i,\bu ^*) \\
 & = & b(\be_i,\bu^n,\bu ^*) + 0 \\
 & \ne & 0
\end{eqnarray*}
in general, and similarly,
\begin{eqnarray*}
 b(\widetilde{\bphi_i},\bu^n,\bu ^*) - b(\bu^n,\widetilde{\bphi_i},\bu ^*)
 & = & b(\bphi_i,\bu^n,\bu ^*) - b(\bu^n,\bphi_i,\bu ^*) \\
 & \ne & 0,
\end{eqnarray*}

Hence we conclude that the linearization \eqref{emacL1}-\eqref{emacL2} conserves energy, but not momentum or angular momentum.

\subsection{The Newton linearization for EMAC}

A common approach to the linearization of flow problems is to use a fixed number (one or two) of the Newton iterations to solve for  solution of a fully implicit scheme.  To make our argument more clear, we consider the Crank-Nicolson scheme and one Newton iteration to resolve the non-linearity. This can be readily seen to be equivalent to the following \textit{linear} time-stepping method:
\begin{align}
\left(\frac{\bu^n-\bu^{n-1}}{\Delta t},\bv\right)  + \nu (\nabla \bu^{n+\frac12},\nabla \bv)
+  2(D(\bu^*)\bu^{n+\frac12}+ D(\bu^{n+\frac12})\bu^*&- D(\bu^*)\bu^*,\bv) \nonumber \\
+ ((\div \bu^*)\bu^{n+\frac12}+ (\div \bu^{n+\frac12})\bu^*-(\div \bu^*)\bu^*,\bv) - (p^{n},\div \bv)
& = (\blf^{n+\frac12},\bv), \label{emacN1} \\
(\div \bu^n,q) & = 0, \label{emacN2}
\end{align}
where $\bu^{n+\frac12}=\frac12(\bu^n+\bu^{n-1})$ and $\bu^*$ is some known approximation of $\bu^{n+\frac12}$.  We choose $\bu^*=\frac32\bu^{n-1} - \frac12\bu^{n-2}$ for a second order method.

To analyze conservation properties of \eqref{emacN1}--\eqref{emacN2}, we need only consider if the nonlinear terms vanish when $\bv=\bu^{n+\frac12},\ \widetilde{\be_i},\ \widetilde{\bphi_i}$.

We start with the momentum conservation and choose $\bv = \widetilde{\be_i}$.  The six nonlinear terms reduce via
\begin{eqnarray*}
&& 2(D(\bu^*)\bu^n+ D(\bu^n)\bu^*-D(\bu^*)\bu^*,\widetilde{\be_i}) + ((\div \bu^*)\bu^n+(\div \bu^n)\bu^*-(\div \bu^n)\bu^*,\widetilde{\be_i}) \\
&& =2(D(\bu^*)\bu^n+D(\bu^n)\bu^*-D(\bu^*)\bu^*,\be_i) + ((\div \bu^*)\bu^n+(\div \bu^n)\bu^*-(\div \bu^*)\bu^*,\be_i) \\
&& =  b(\bu^*,\bu^n,\be_i) + b(\be_i,\bu^n,\bu^*)   + b(\bu^n,\bu^*,\be_i) + b(\be_i,\bu^*,\bu^n)    + b(\bu^*,\bu^*,\be_i) + b(\be_i,\bu^*,\bu^*)   \\
&& \hspace{.5in} + ((\div \bu^*)\bu^n,\be_i) + ((\div \bu^n)\bu^*,\be_i) - ((\div \bu^*)\bu^*,\be_i)
\end{eqnarray*}
Notice that since $\be_i$ is constant, identity \eqref{vecid1} gives $b(\bv,\bw,\be_i)  = - ((\div \bv)\bw,\be_i)$ for any $\bv,\bw\in \bX$, and thus 6 terms of the 9 term expansion drop, leaving just 3 terms,
\begin{eqnarray*}
b(\be_i,\bu^n,\bu^*)   + b(\be_i,\bu^*,\bu^n)    + b(\be_i,\bu^*,\bu^*)  = 0
\end{eqnarray*}
with the last step following because $\be_i$ is constant and thus divergence-free, we have again from \eqref{vecid1} that $b(\be_i,\bu^*,\bu^*) =0$ and $b(\be_i,\bu^n,\bu^*) = - b(\be_i,\bu^*,\bu^n) $.\\

To see angular momentum conservation, choose $\bv=\widetilde{\bphi_i}$.  Similar to momentum, reducing and expanding the nonlinear terms gives
\begin{eqnarray*}
&& 2(D(\bu^*)\bu^n+D(\bu^n)\bu^*-D(\bu^*)\bu^*,\bphi_i) + ((\div \bu^*)\bu^n+(\div \bu^n)\bu^*-(\div \bu^*)\bu^*,\bphi_i) \\
& = & b(\bu^*,\bu^n,\bphi_i) + b(\bphi_i,\bu^n,\bu^*)   + b(\bu^n,\bu^*,\bphi_i) + b(\bphi_i,\bu^*,\bu^n)    + b(\bu^*,\bu^*,\bphi_i) + b(\bphi_i,\bu^*,\bu^*)   \\
&& \hspace{.5in} + ((\div \bu^*)\bu^n,\bphi_i) + ((\div \bu^n)\bu^*,\bphi_i) - ((\div \bu^*)\bu^*,\bphi_i).
\end{eqnarray*}
Since $\div \bphi_i =0$, we have from \eqref{vecid1} that
\[
 b(\bphi_i,\bu^n,\bu^*)   ) + b(\bphi_i,\bu^*,\bu^n)    + b(\bphi_i,\bu^*,\bu^*) = 0.
\]
For the remaining 6 terms, consider them 2 at a time.  For the two that have $\bu^*$ in their first argument,
\begin{eqnarray*}
b(\bu^*,\bu^n,\bphi_i) + ((\div \bu^*)\bu^n,\bphi_i) & = & -b(\bu^*,\bphi_i,\bu^n)  - ((\div \bu^*)\bu^n,\bphi_i) + ((\div \bu^*)\bu^n,\bphi_i) \\
& = &  -b(\bu^*,\bphi_i,\bu^n).
\end{eqnarray*}
However, since  $(\nabla \bphi_i)$ is a skew symmetric matrix, we can write
\begin{eqnarray*}
-b(\bu^*,\bphi_i,\bu^n)
& = & -((\nabla \bphi_i)\bu^*,\bu^n) \\
& = & -(\bu^*,(\nabla \bphi_i)^T\bu^n) \\
& = & (\bu^*,(\nabla \bphi_i)\bu^n) \\
& = & b(\bu^n,\bphi_i,\bu^*),
\end{eqnarray*}
and thus
\[
b(\bu^*,\bu^n,\bphi_i) + ((\div \bu^*)\bu^n,\bphi_i) = -b(\bu^*,\bphi_i,\bu^n) = 0.
\]
Similarly, we can show for the other four terms in the expansion above that
\[
b(\bu^n,\bu^*,\bphi_i) + ((\div \bu^n)\bu^*,\bphi_i) =0, \ \ \ b(\bu^*,\bu^*,\bphi_i) + ((\div \bu^*)\bu^*,\bphi_i)=0.
\]
Hence we have found that the Newton linearization for EMAC \eqref{emacN1}-\eqref{emacN2} conserves momentum and angular momentum.

Now we turn to study the energy balance. Although energy is not exactly conserved unless the Newton iteration converges, we shall see that under assumption that  the discrete solution behaves sufficiently regular, the error introduced to the energy balance on each time step is of higher order compared to the approximation error (for the scheme in \eqref{emacN1}--\eqref{emacN2} it is $O(|\Delta t|^5)$ vs  $O(|\Delta t|^3)$) at a given time step.
To see this, we set $\bv=\bu^{n+\frac12}$, which gives for the nonlinear terms
\begin{multline}
N\,:=\,2(D(\bu^*)\bu^{n+\frac12}+D(\bu^{n+\frac12})\bu^*-D(\bu^*)\bu^*,\bu^{n+\frac12})\\ + ((\div \bu^*)\bu^{n+\frac12}+(\div \bu^{n+\frac12})\bu^*-(\div \bu^*)\bu^*,\bu^{n+\frac12}). \label{NE1}
\end{multline}
From \eqref{econs}, we know that for $\bw\in \bX$, $2(D(\bw)\bw,\bw) + ((\div \bw)\bw,\bw)=0$. Setting $\bw=\bu^{n+\frac12}$ and subtracting from  \eqref{NE1} gives after re-grouping,
\[
N=2(D(\bu^*-\bu^{n+\frac12})(\bu^{n+\frac12}-\bu^*),\bu^{n+\frac12}) + ((\div (\bu^*-\bu^{n+\frac12}))(\bu^{n+\frac12}-\bu^*),\bu^{n+\frac12}).
\]
Now we substitute  $\bu^*=\frac32\bu^{n-1} - \frac12\bu^{n-2}$ to find
\begin{equation}\label{Aux1}
N=-\frac{|\Delta t|^4}2(D(\left[\bu\right]^{n-1}_{tt})\left[\bu\right]^{n-1}_{tt},\bu^{n+\frac12}) - \frac{|\Delta t|^4}4((\div (\left[\bu\right]^{n-1}_{tt})(\left[\bu\right]^{n-1}_{tt})),\bu^{n+\frac12}),
\end{equation}
where  $\left[\bu\right]^{n-1}_{tt}=\frac{\bu^n-2\bu^{n-1}+\bu^{n-2}}{(\Delta t)^2}$.
If we \textit{assume} that the second finite differences $\left[\bu\right]^{n-1}_{tt}$ are bounded in suitable norms (which is the case if the numerical solution approximate a smooth true NSE solution in those norms
with $O(|\Delta t|^2)$ order), then we find that $N$ is $O(|\Delta t|^4)$--small,
\[
|N|\le C(\bu)|\Delta t|^4.
\]
From this and \eqref{emacN1}--\eqref{emacN2} with $\bv=\bu^{n+\frac12}$, we obtain the energy balance of the
following form
\[
E_{\rm kin}(t^n)+\Delta t\nu\|\bu^{n+\frac12}\|^2 =E_{\rm kin}(t^{n-1})+\Delta t(\blf^{n+\frac12},\bu^{n+\frac12})+O(|\Delta t|^5).
\]
Thus, on every time step an error of the order $O(|\Delta t|^5)$ is introduced to the energy balance.
We should note again that the crucial assumption was that the solution of the discrete problem is sufficiently regular, e.g. approximates a smooth solution. Proving convergence or regularity property of solutions to linearized EMAC schemes is an open research topic.

\section{Numerical Experiments}

We now present results for several numerical test problems, which compare the various EMAC schemes proposed above, and also demonstrate the ability of EMAC to give good results on difficult test problems.

\subsection{Gresho Problem}

\begin{figure}[!h]
\begin{center}
\includegraphics[width=.31\textwidth,height=0.31\textwidth, viewport=70 45 550 400, clip]{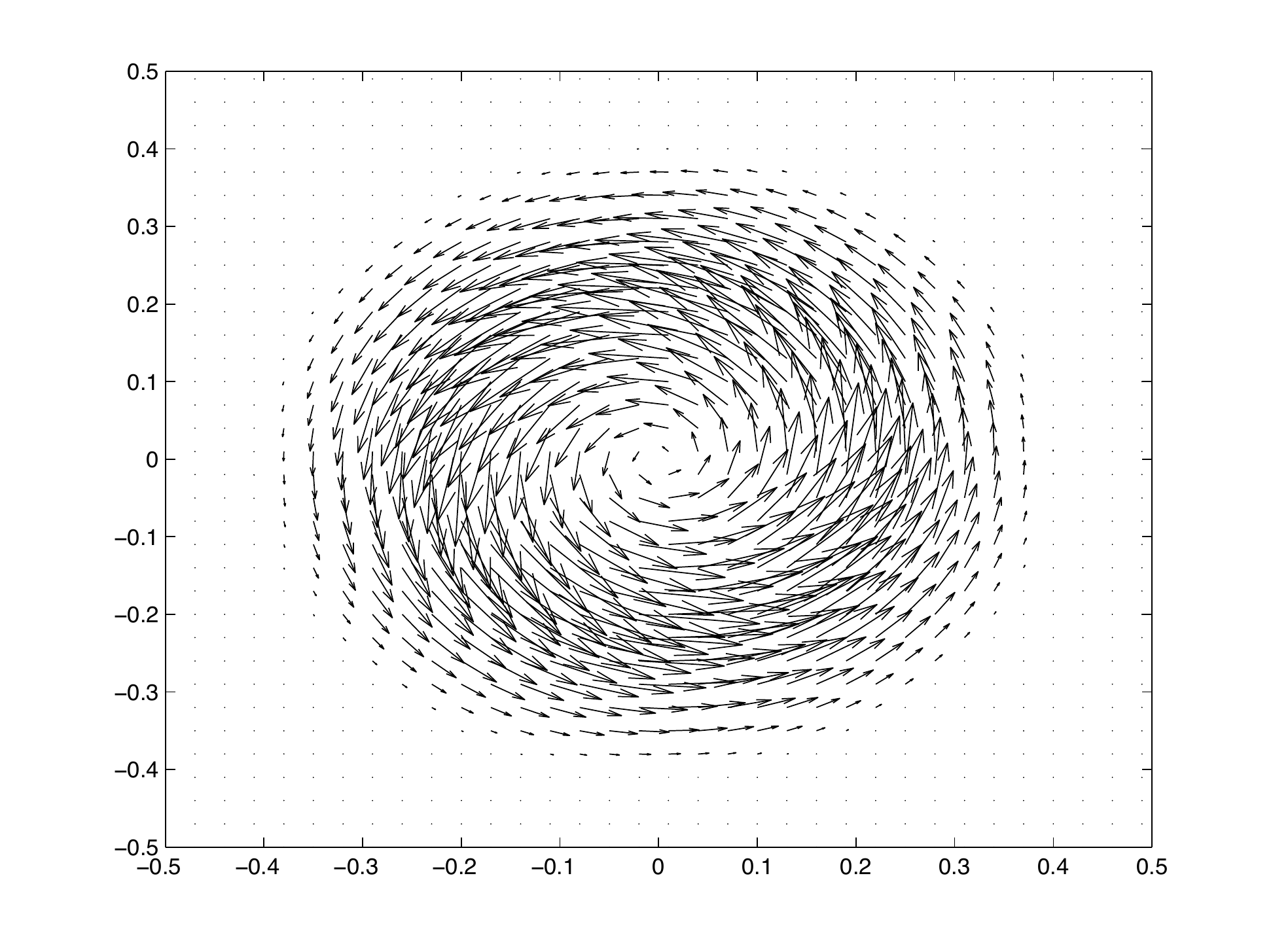}
\includegraphics[width=.35\textwidth,height=0.34\textwidth, viewport=70 45 550 440, clip]{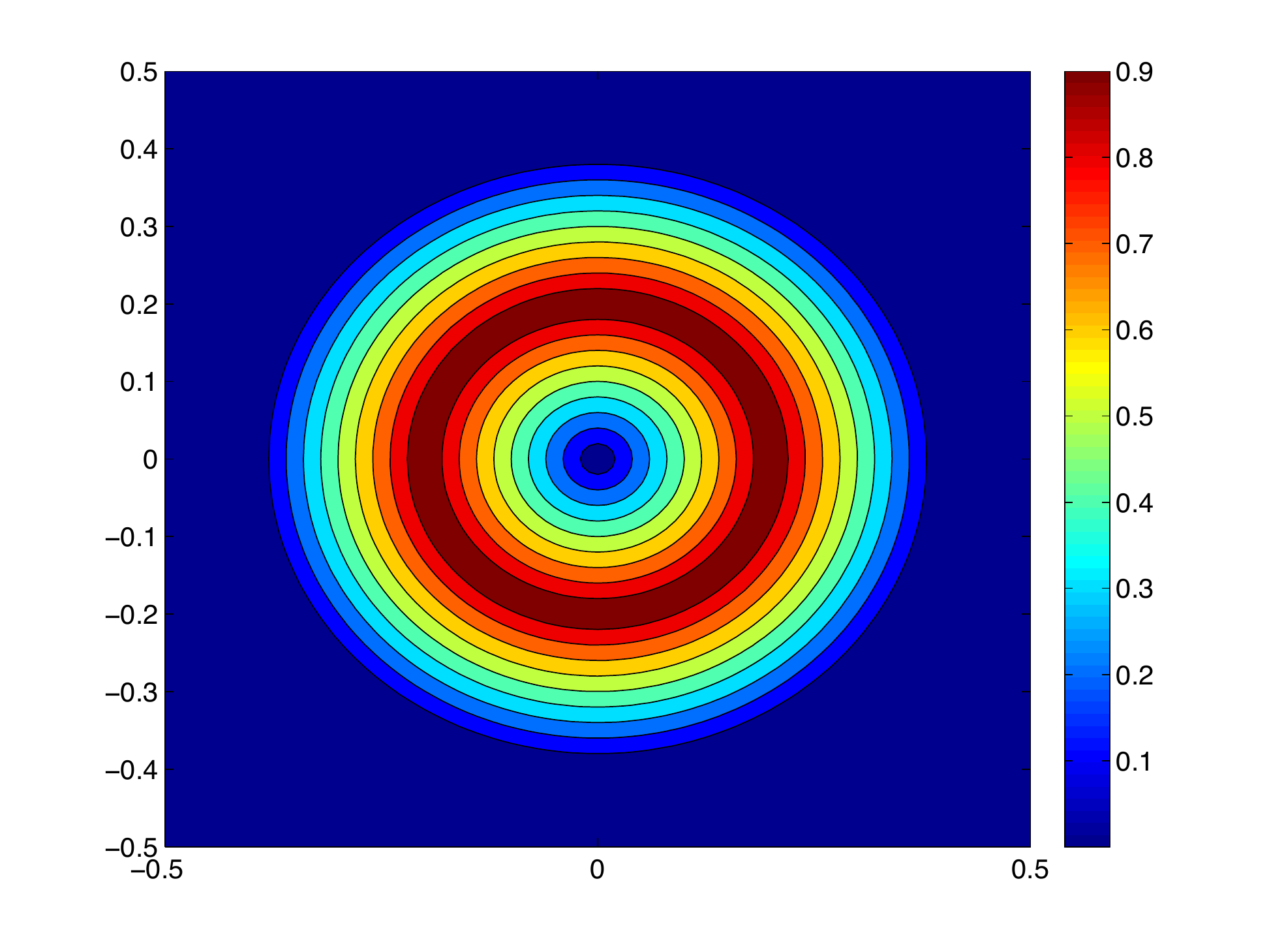}
\end{center}
\caption{\label{vortexmesh}
Shown above is the velocity solution for the Gresho problem as a vector plot (left) and speed contour plot (right). }
\end{figure}

The first test problem we consider is the Gresho problem, which is also known as the `standing vortex problem' \cite{TMRS92,LW03,G90}.  The problem starts with an initial condition that is an exact solution of the steady Euler equations, and then run with a time stepping method using $\blf = {\bf 0}$, $\nu=0$, and no penetration boundary conditions.  A correct solution will not change in time.

On the domain $\Omega=(-.5,.5)^2$, with $r=\sqrt{x^2+y^2}$, velocity and pressure that are solutions of the steady Euler equations are given by
\[
\begin{aligned}
r \le 0.2~:& \left\{\begin{split} \bu &= \left( \begin{array}{c} -5y \\ 5x \end{array} \right) \\
 p&=12.5r^2 + C_1 \end{split}\right.,\\
r>0.4 ~:&  \left\{ \begin{split} \bu& = \left( \begin{array}{c} 0 \\ 0 \end{array} \right)\\ p&=0\end{split}\right.,  \\
0.2 \le r \le 0.4~:& \left\{\begin{split}\bu &= \left( \begin{array}{c}  \frac{-2y}{r}+5y \\ \frac{2x}{r} - 5x \end{array} \right)   \\
 p&=12.5r^2 - 20r + 4\log(r) + C_2 \end{split}\right. ,
\end{aligned}
\]
where
\[
C_2 = (-12.5)(0.4)^2 + 20(0.4)^2 - 4\log(0.4),\ C_1 = C_2 - 20(0.2) + 4\log(0.2).
\]
Vector and speed plots are given for the velocity in figure \ref{vortexmesh}.

We choose this Euler velocity solution $\bu$ to be the initial condition for our numerical simulations, and note that an accurate scheme will preserve the solution in time.  Moreover, it is also a good test for conservation properties, as there is no viscosity or external forcing, and the boundaries do not influence the solution (unless significant error is already present).

We compute solutions to the Gresho problem using the different EMAC formulations, together with Crank-Nicolson time stepping, with ${\bf f}={\bf 0}$, $\nu=0$, and no-penetration boundary conditions up to T=10.  We computed using $(P_2,P_1)$ Taylor-Hood element on a 48x48 uniform mesh, with a time step of $\Delta t=0.01$.  The formulations we used were fully nonlinear (using Newton to resolve the the nonlinear problem at each time step [$\bH^1$ norm of step size $<10^{-8}$]), the energy stable skew-symmetric linearization, and the Newton linearization with 1, 2, and 3 steps of Newton at each time step (it typically took the nonlinear EMAC scheme 3 or 4 Newton iterations per time step to converge).

\begin{figure}[!h]
\begin{center}
\includegraphics[width=.49\textwidth,height=0.22\textwidth, viewport=60 0 900 320, clip]{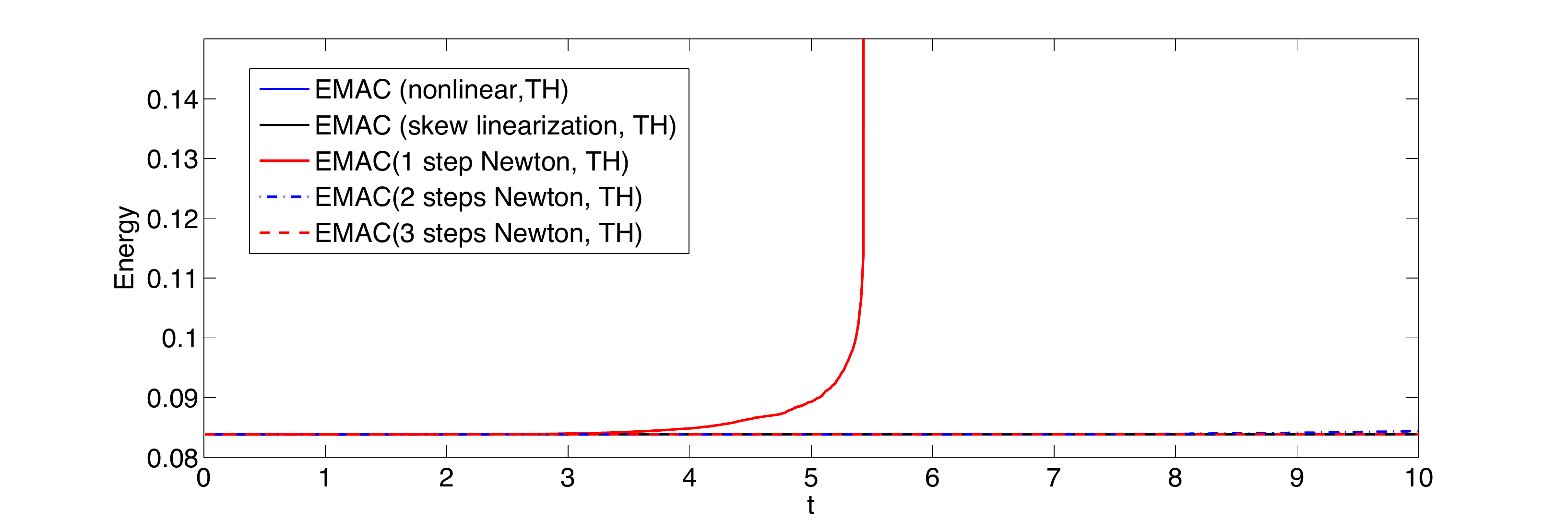}
\includegraphics[width=.49\textwidth,height=0.22\textwidth, viewport=60 0 900 320, clip]{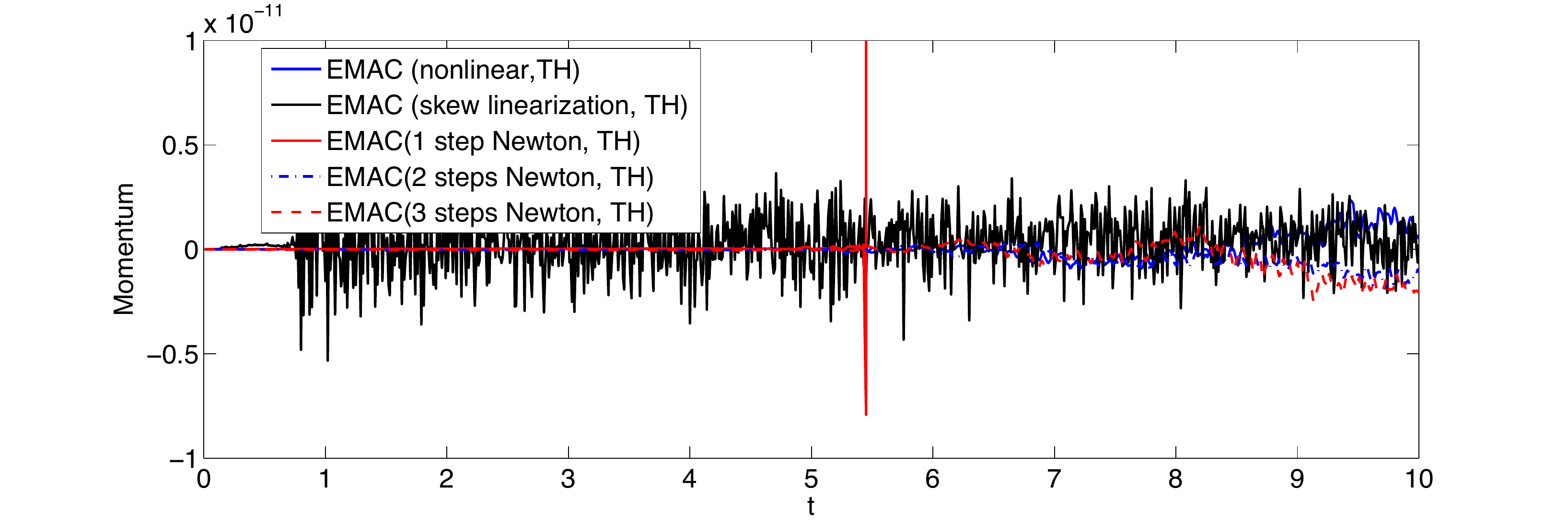}
\includegraphics[width=.49\textwidth,height=0.22\textwidth, viewport=60 0 900 320, clip]{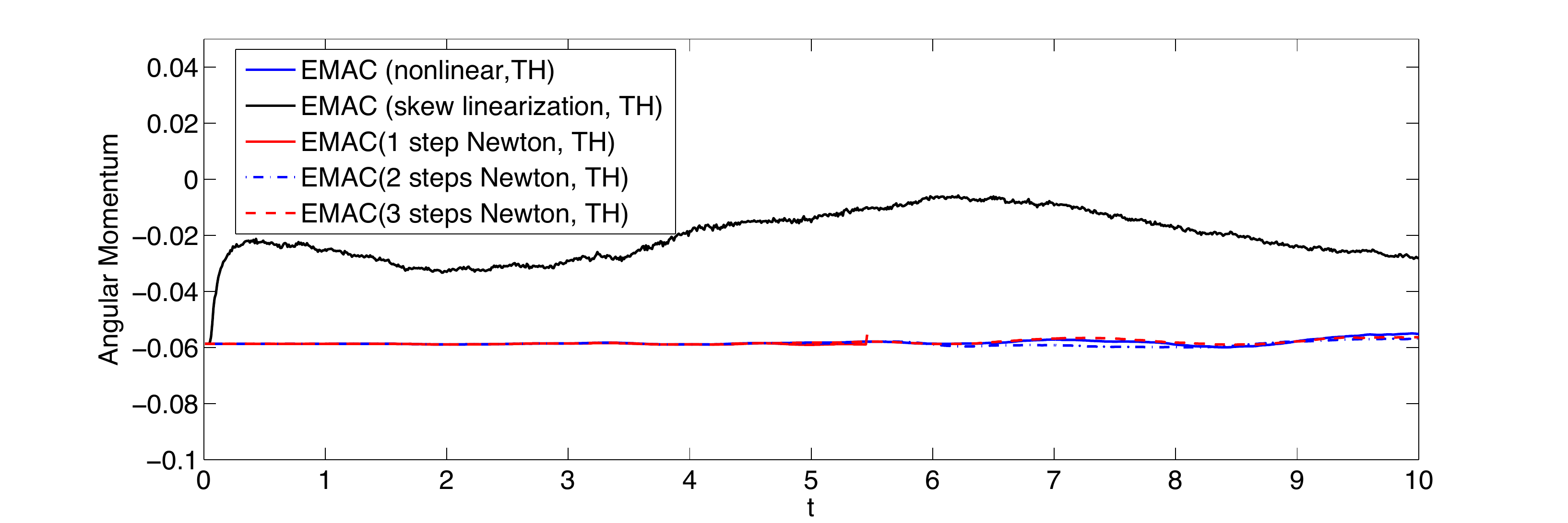}
\includegraphics[width=.49\textwidth,height=0.22\textwidth, viewport=60 0 900 320, clip]{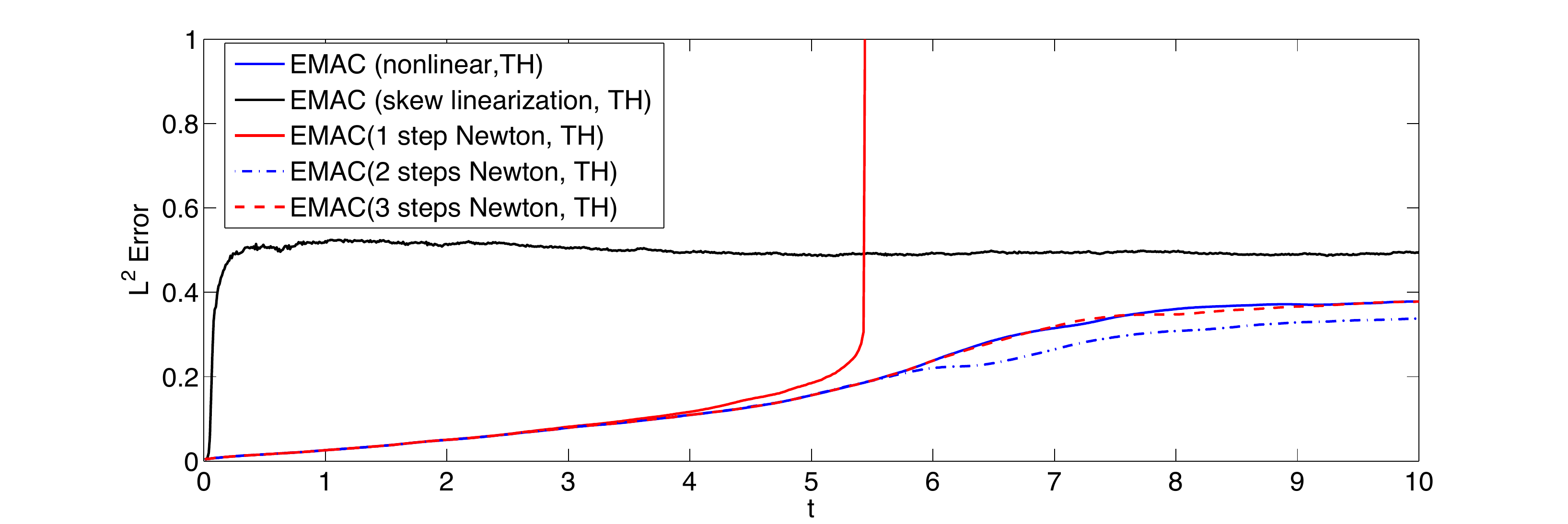}
\end{center}
\caption{\label{Gresho5}
Shown above are plots of time versus energy, momentum, angular momentum, and $L^2(\Omega)$ velocity error, for the various EMAC linearizations in the Gresho problem, using Taylor-Hood elements.}
\end{figure}

Results for the Gresho tests are shown in figure \ref{Gresho5}, as plots of energy, momentum, angular momentum, and $L^2(\Omega)$ error versus time.  We observe all the schemes conserve energy, except for the one that uses just 1 step of Newton at each time step, which is consistent with our analysis above that shows the Newton method conserves energy only up to the level that the Newton method converged.  Hence the 1-step-Newton scheme was not energy stable, and blew up around t=5.
For momentum, we see that only the 1-step-Newton method did not conserve momentum well.  Even though we show above that the skew linearization need not conserve momentum, it is aided in this case by the conservation of mass enforcement along with the zero boundary condition, and is able to attain good momentum conservation.  For angular momentum conservation, the unstable method, of course, did not conserve angular momentum (not plotted after blow up), nor did the skew-symmetric linearization's solution.  The angular momentum of the skew symmetric linearization became poor almost immediately.  The angular momentum of the other schemes (2,3, or more steps of Newton) all conserved angular momentum quite well.  Finally, for $L^2(\Omega)$ error, the schemes that conserved angular momentum had much better error, that the skew linearization.  The overall accuracy of the EMAC Newton linearizations were about the same, if 2 or more Newton steps were taken.

\begin{figure}[!h]
\begin{center}
\includegraphics[width=.49\textwidth,height=0.22\textwidth, viewport=60 0 900 320, clip]{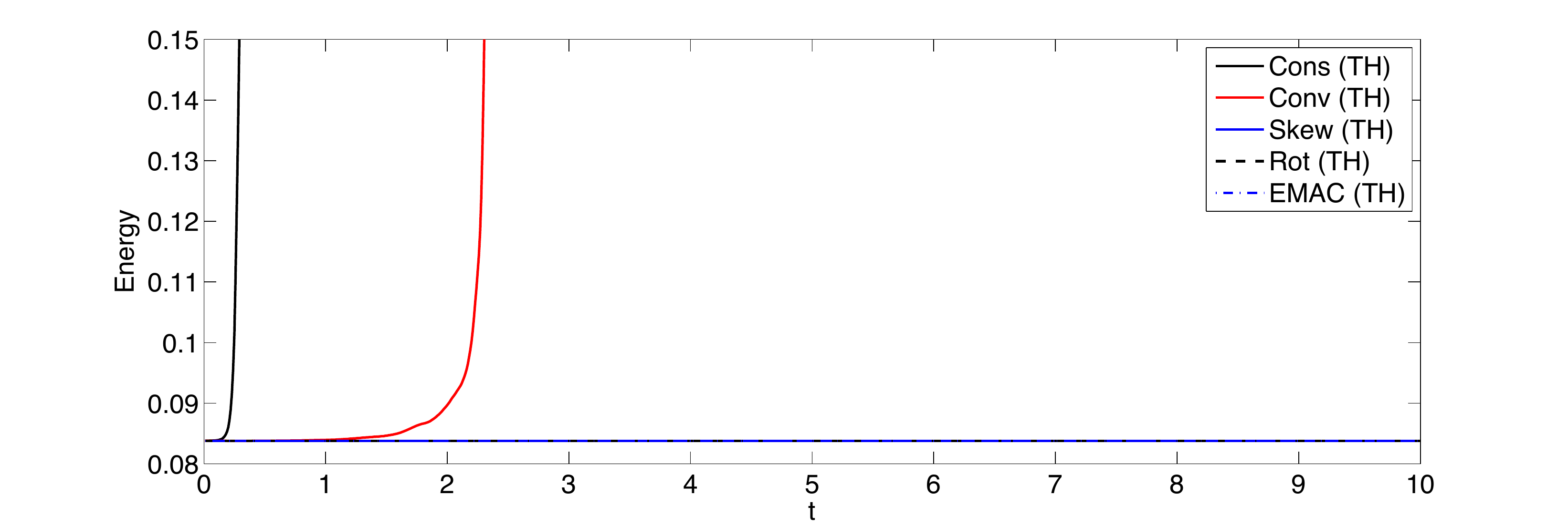}
\includegraphics[width=.49\textwidth,height=0.22\textwidth, viewport=60 0 900 320, clip]{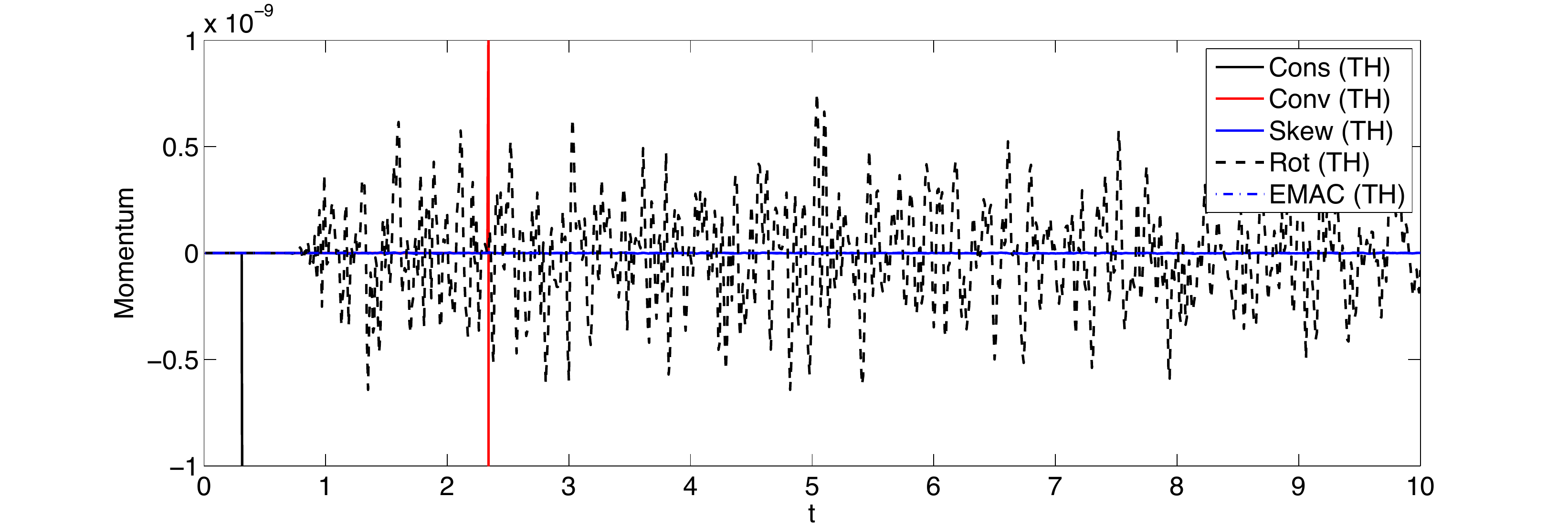}
\includegraphics[width=.49\textwidth,height=0.22\textwidth, viewport=60 0 900 320, clip]{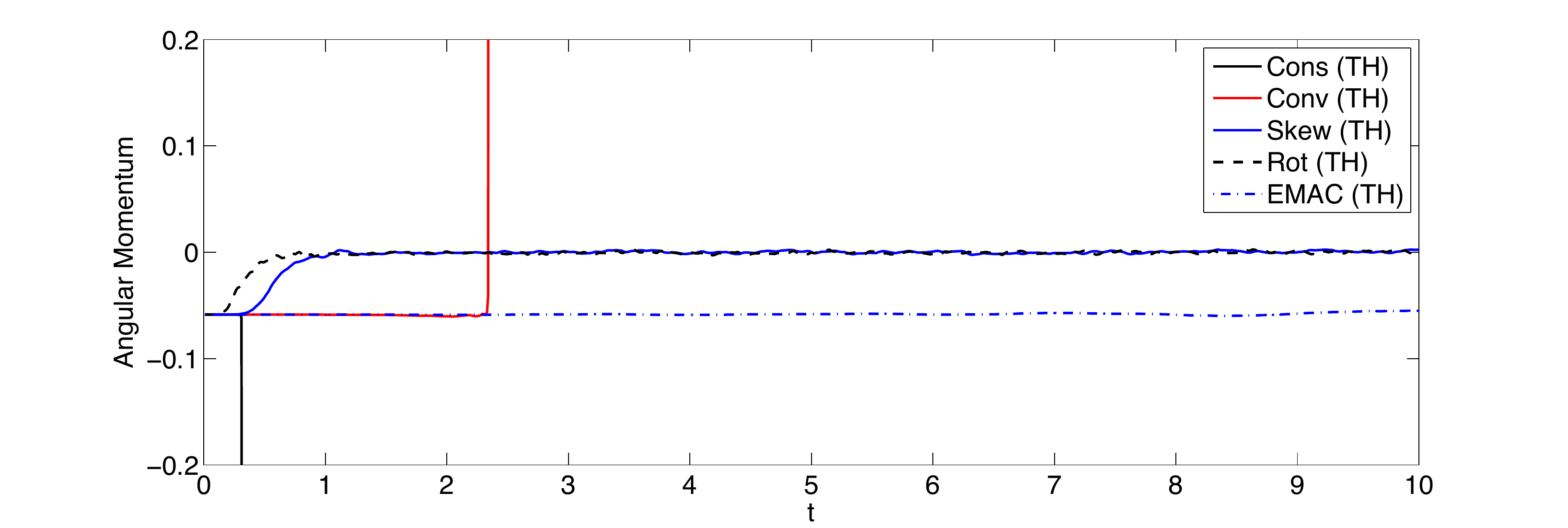}
\includegraphics[width=.49\textwidth,height=0.22\textwidth, viewport=60 0 900 320, clip]{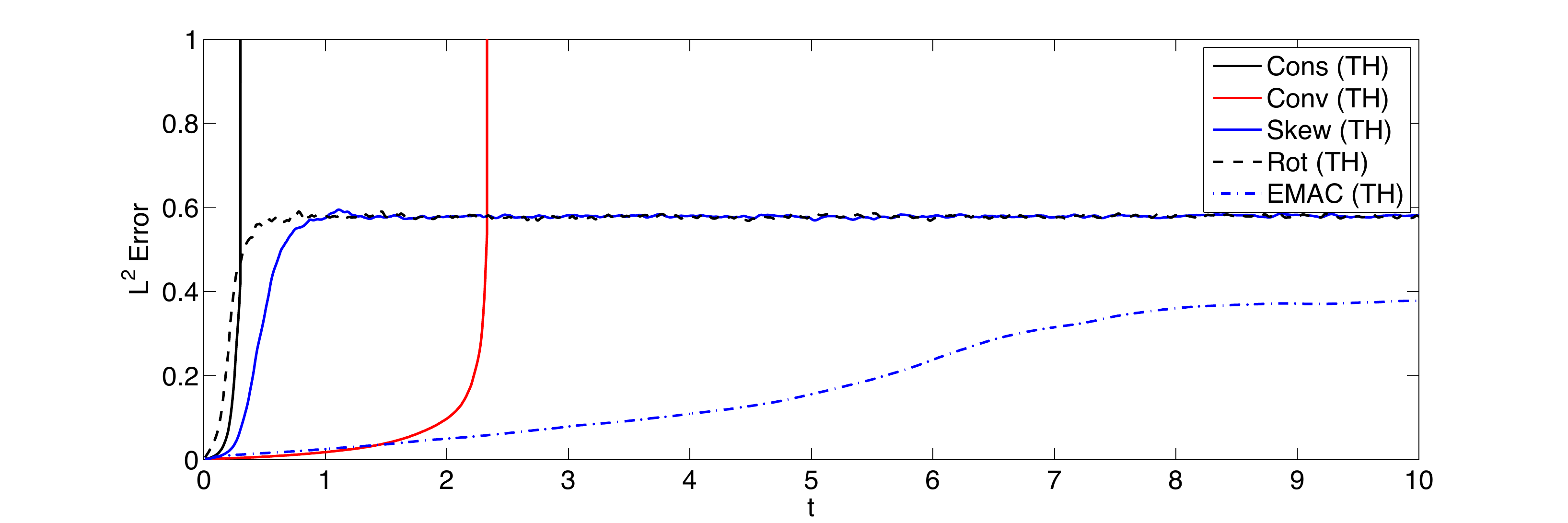}
\end{center}
\caption{\label{Gresho2}
Shown above are plots of time versus energy, momentum, angular momentum, and $L^2(\Omega)$ velocity error, for the various formulations in the Gresho problem, using Taylor-Hood elements.}
\end{figure}

For an additional comparison, we also computed solutions for the other common formulations of the nonlinearity, using analogous numerical schemes (same discretization parameters).  The other formulations we considered were convective (CONV), skew-symmetric (SKEW), conservative (CONS), rotational (ROT), which are computed with the same scheme as \eqref{emac1}-\eqref{emac2}, but with nonlinear terms replaced by:
\begin{align*}
CONV: \  & (\bu^n \cdot \nabla \bu^n,\bv) \\
SKEW:\  & (\bu^n \cdot \nabla \bu^n,\bv) + \frac12 ((\div \bu^n)\bu^n,\bv)\\
CONS:\  & (\bu^n \cdot \nabla \bu^n,\bv) +  ((\div \bu^n)\bu^n,\bv)\\
ROT:\  & ((\nabla \times \bu^n)\times \bu^n,\bv)
\end{align*}
For each of these schemes, the nonlinear problem is solved at each time step, with the same tolerance as for the EMAC nonlinear scheme.  Results are shown in figure \ref{Gresho2}, along with results for the (nonlinear) EMAC scheme.  We observe that EMAC outperforms all the other formulations, especially over longer times (i.e. CONV is competitive until it becomes unstable).  Since the 2 step Newton linearization of EMAC gave about the same results as nonlinear EMAC, we can further conclude that the 2 step Newton linearization of EMAC performs much better than each of these other formulations (which solve the nonlinear problem).

\subsection{Lattice Vortex problem}

\begin{figure}[!h]
\begin{center}
\includegraphics[width=.31\textwidth,height=0.25\textwidth, viewport=70 45 550 390, clip]{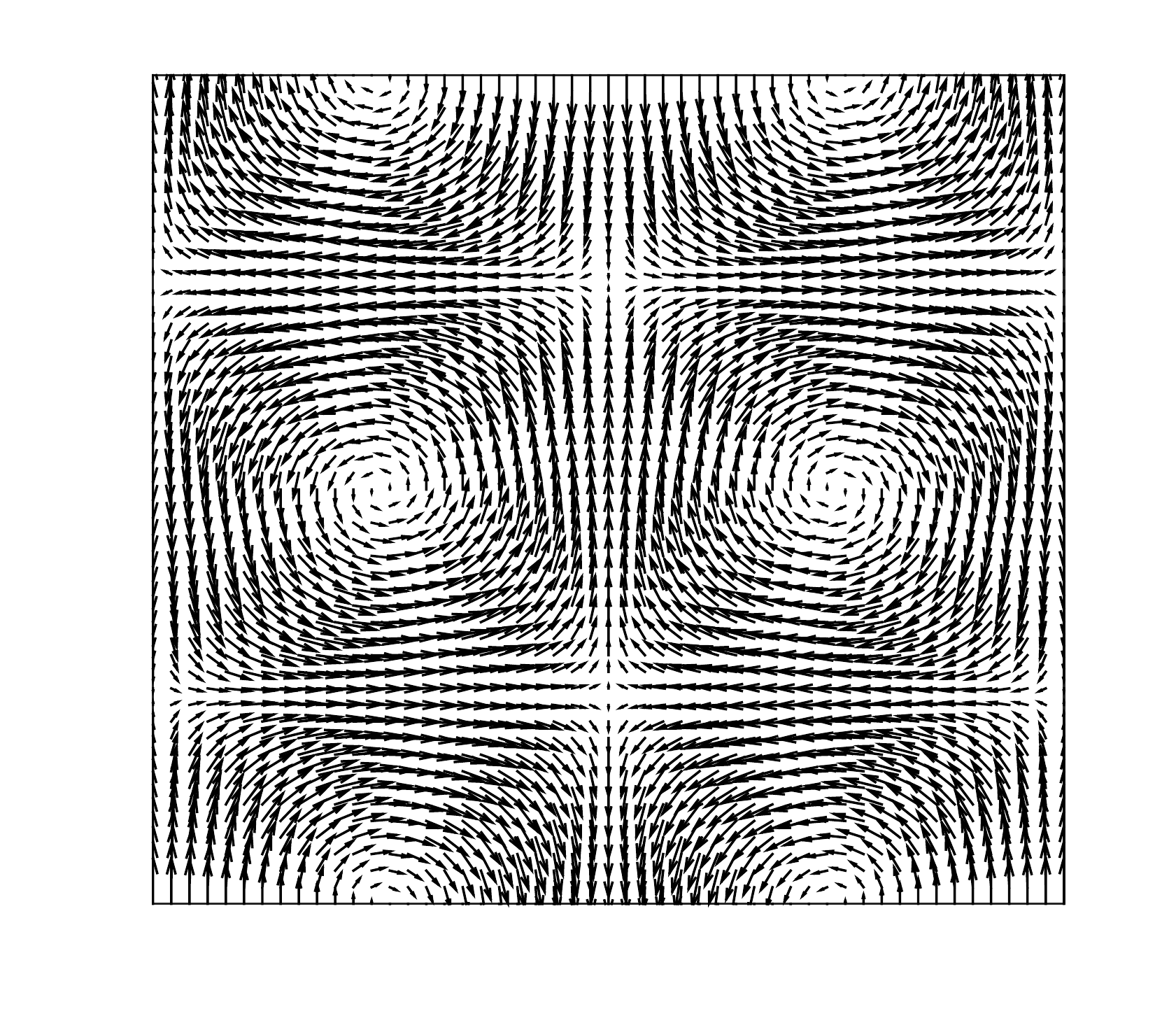}
\includegraphics[width=.31\textwidth,height=0.25\textwidth, viewport=70 45 550 390, clip]{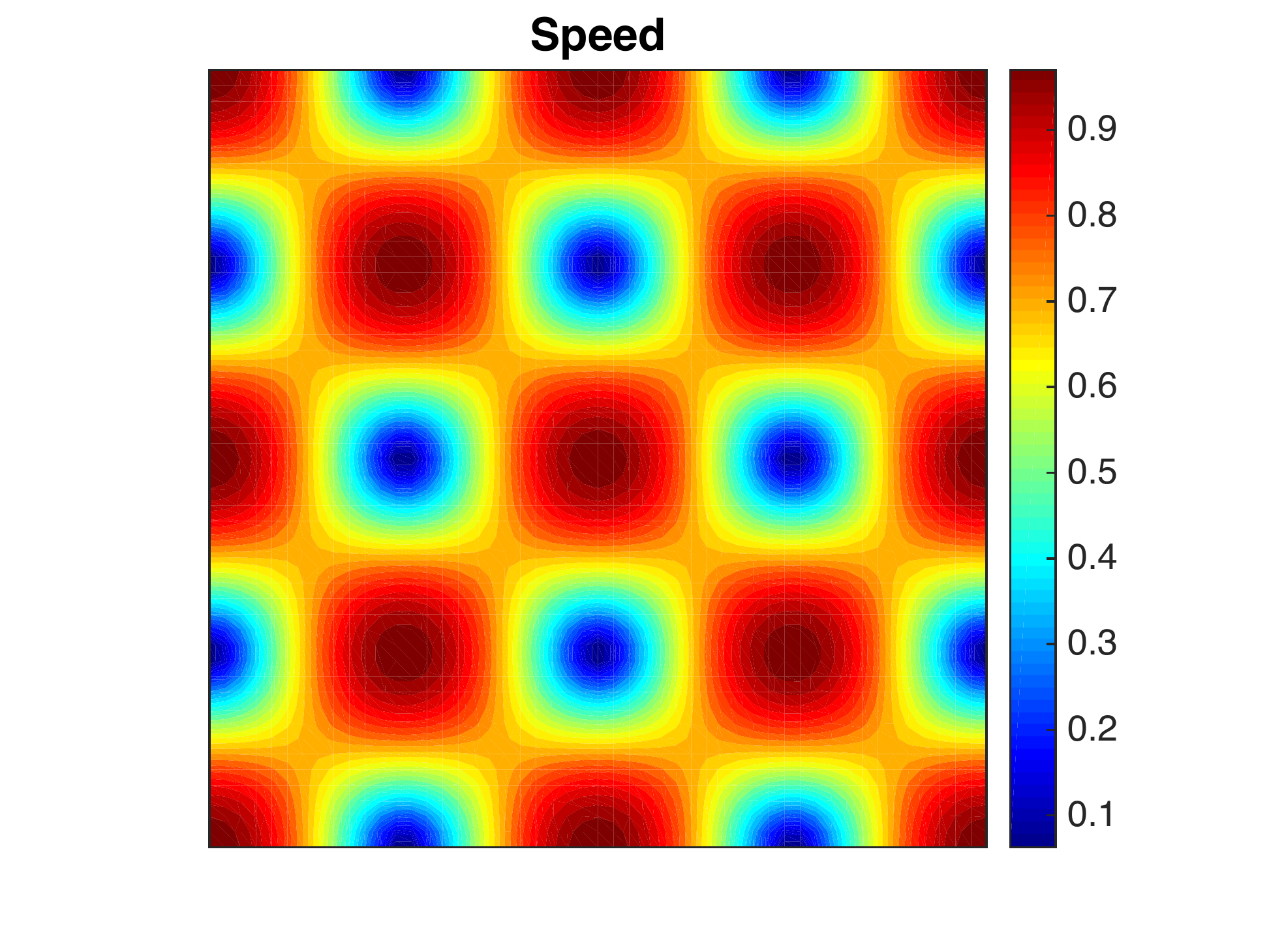}
\includegraphics[width=.31\textwidth,height=0.25\textwidth, viewport=70 45 550 390, clip]{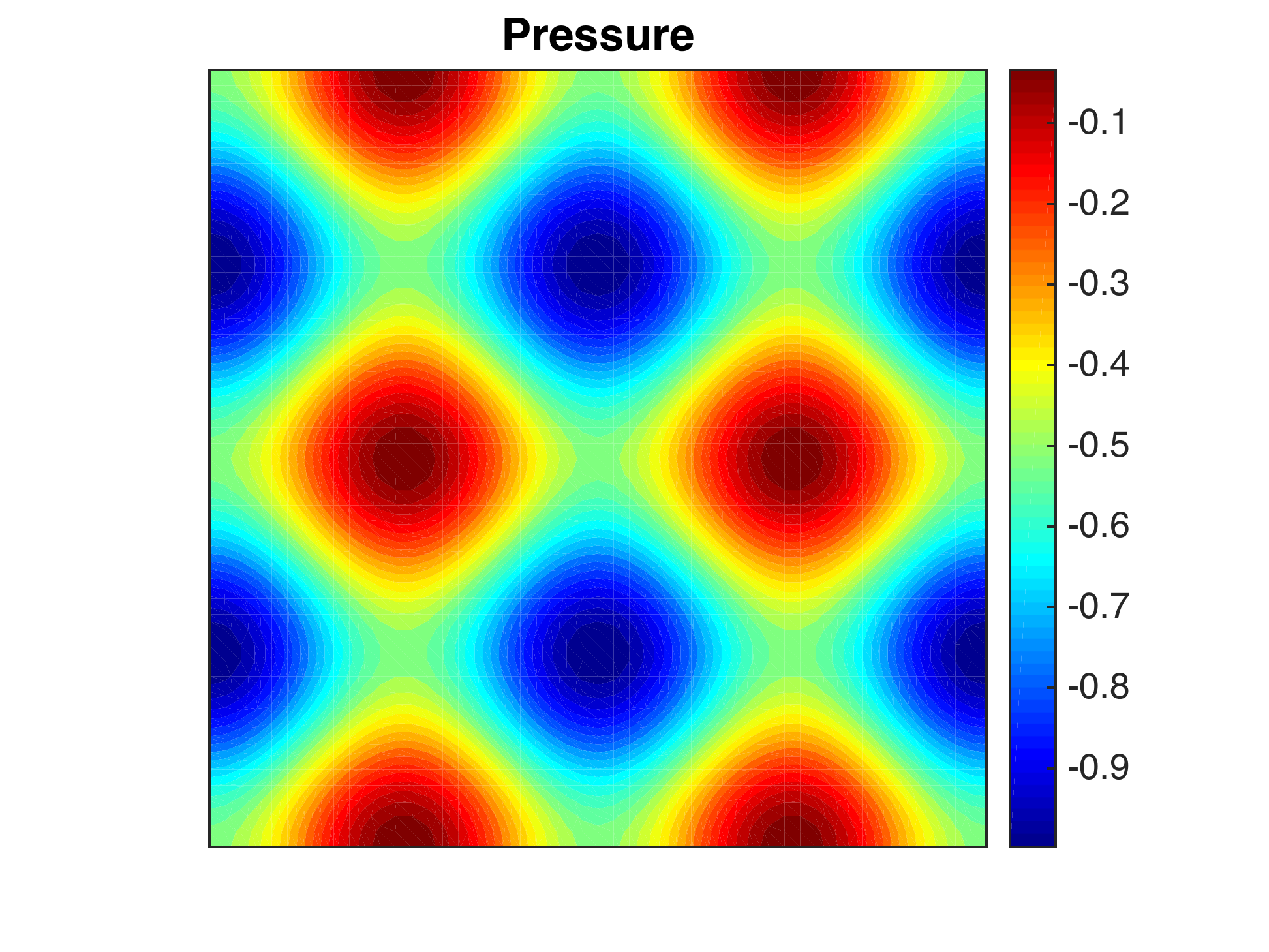}
\end{center}
\caption{\label{lattice}
Shown above is the solution for the lattice vortex problem as a vector plot of velocity (left), speed contour plot (center), and
pressure contour plot (right). }
\end{figure}

The second test problem we consider is the lattice vortex problem from \cite{MB2002,SL17}, with small viscosity.  The true solution takes the form
\[
\bu = \bv e^{-8\nu \pi^2 t}, \ \ p=q e^{-16 \pi^2 t},
\]
where
\[
\bv = \langle \sin (2\pi x) \sin(2\pi y),\ \cos(2\pi x)\cos(2\pi y) \rangle,\ \ \ q = -\frac12 \left( \sin^2(2 \pi x) + \cos^2(2 \pi y) \right).
\]
Note that $(\bu,p)$ is an exact NSE solution with $\blf={\bf 0}$, and $(\bv,q)$ is an exact Euler solution with $\blf={\bf 0}$.  We consider the test problem on $\Omega=(0,1)^2$ with time up to $t=5$, taking $\bu_0=\bu(0)$ as the initial condition, $\blf={\bf 0}$, and $\nu=10^{-7}$.  This is a harder problem than the Gresho problem because there are several spinning vortices whose edges touch, which can be difficult to resolve numerically.  Due to the small viscosity, the true solution will be essentially the same at t=10, although very slightly decayed.  A plot of the initial condition is shown in figure \ref{lattice}.

We simulate this problem using nonlinear schemes for EMAC, CONV, SKEW, CONS, and ROT, as well as the 1, 2, and 3 step Newton linearizations for EMAC.  We used a 1/32 triangular mesh and time step $\Delta t=0.01$, together with Taylor-Hood elements and Crank-Nicolson time stepping.  We strongly enforce Dirichlet boundary conditions for all schemes to be the true velocity solution at the boundary nodes.

\begin{figure}[!h]
\begin{center}
\includegraphics[width=.49\textwidth,height=0.28\textwidth, viewport=0 0 550 420, clip]{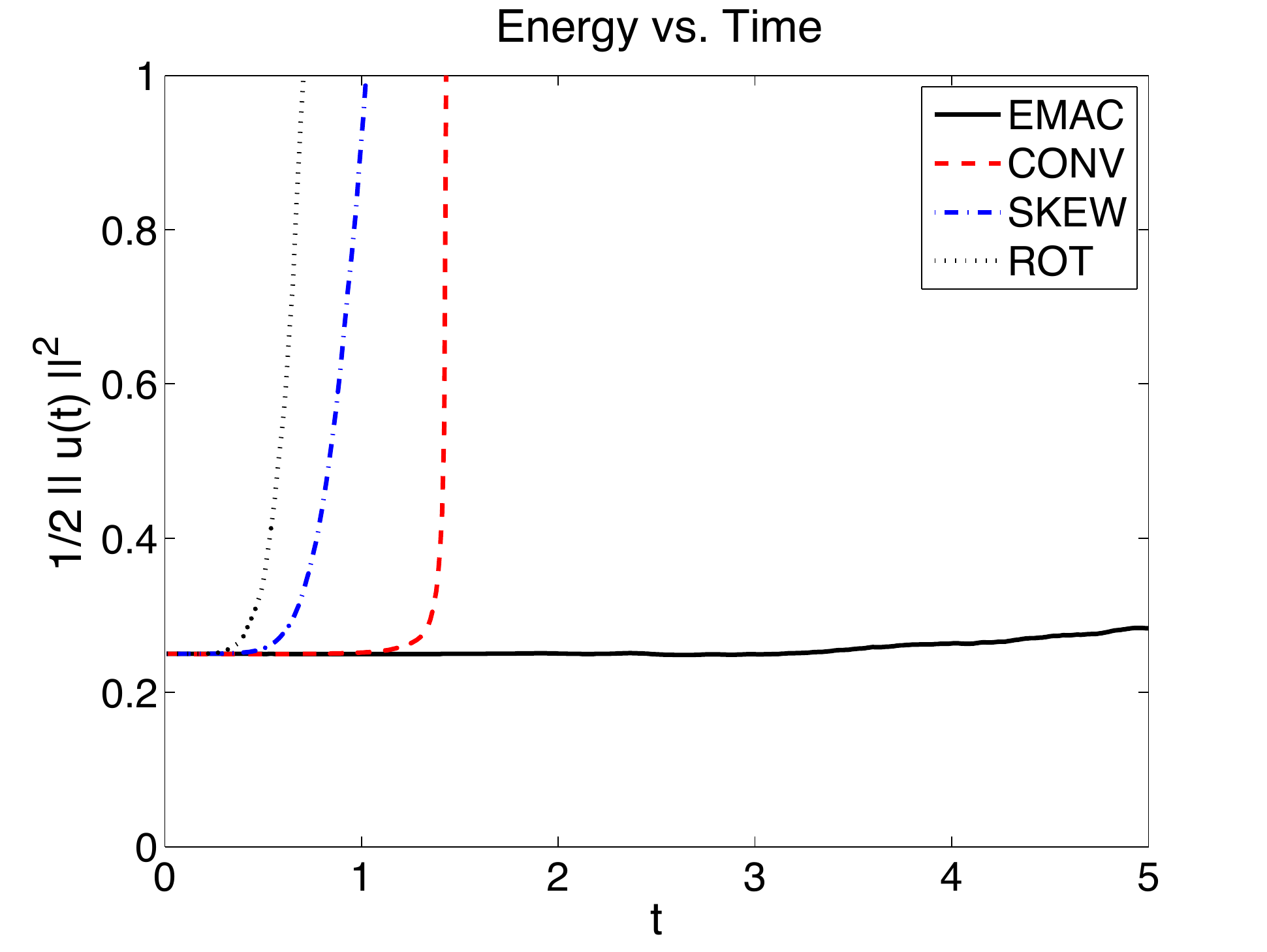}
\includegraphics[width=.49\textwidth,height=0.28\textwidth, viewport=0 0 550 420, clip]{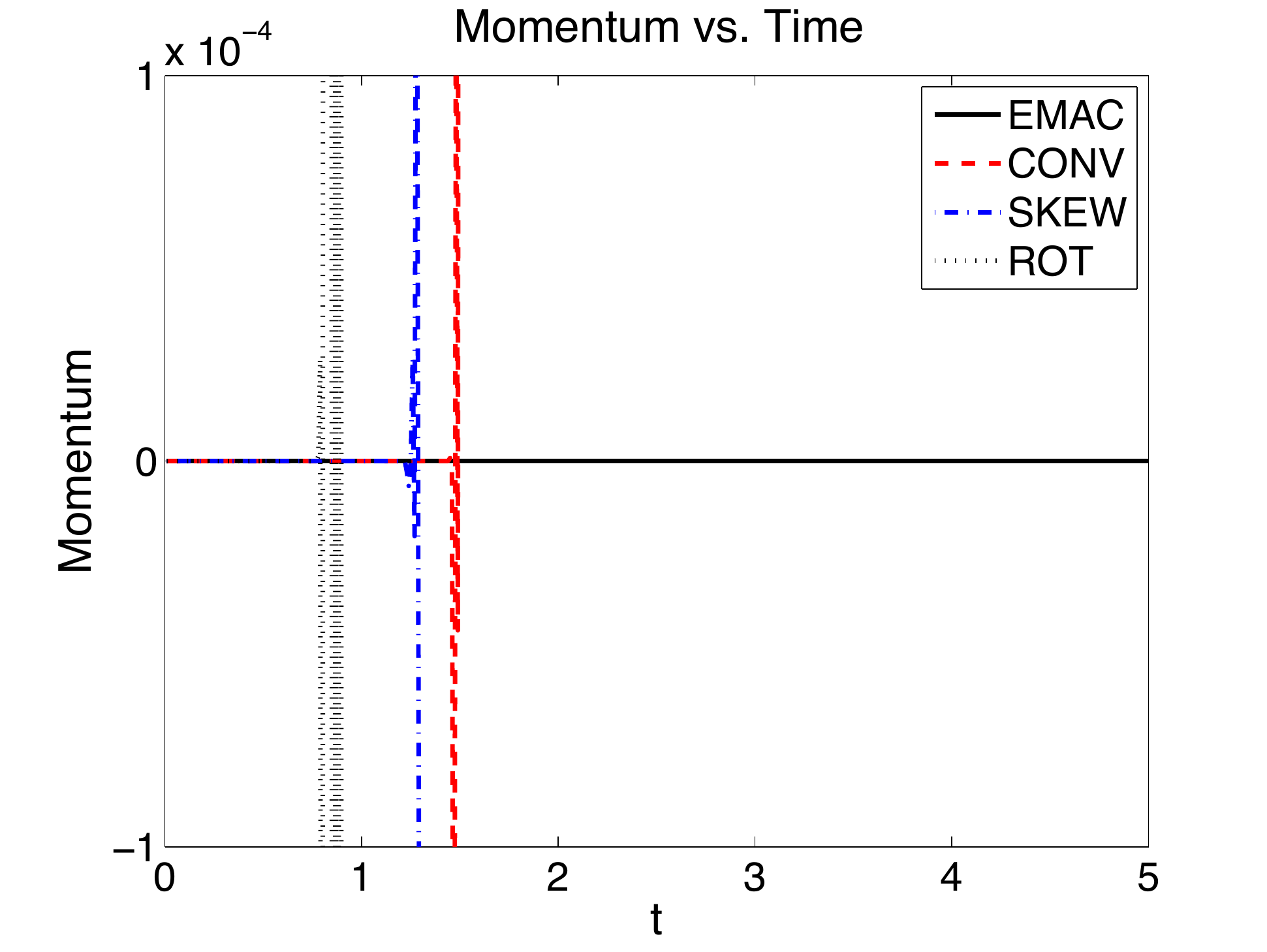}
\includegraphics[width=.49\textwidth,height=0.28\textwidth, viewport=0 0 550 420, clip]{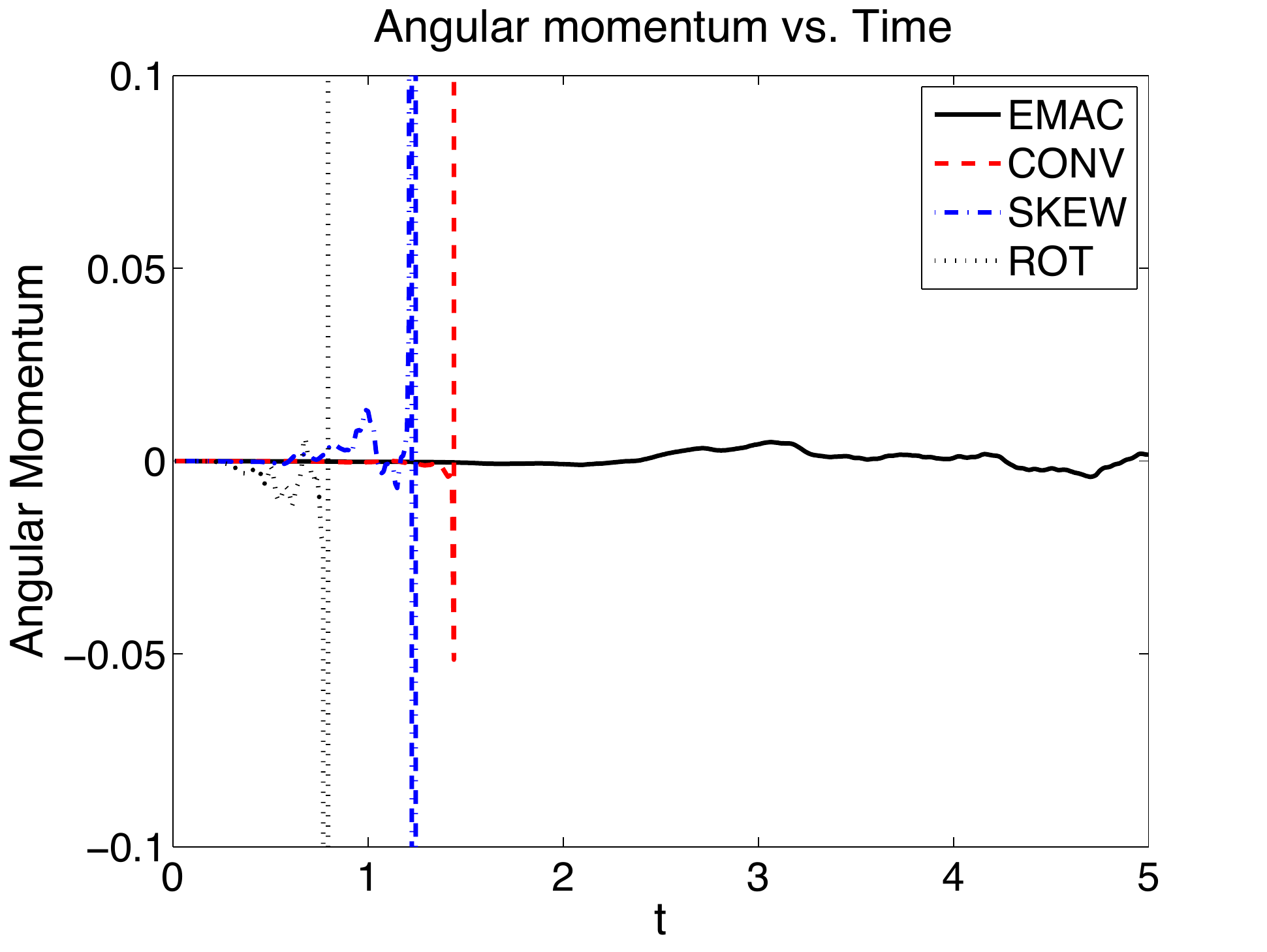}
\includegraphics[width=.49\textwidth,height=0.28\textwidth, viewport=0 0 550 420, clip]{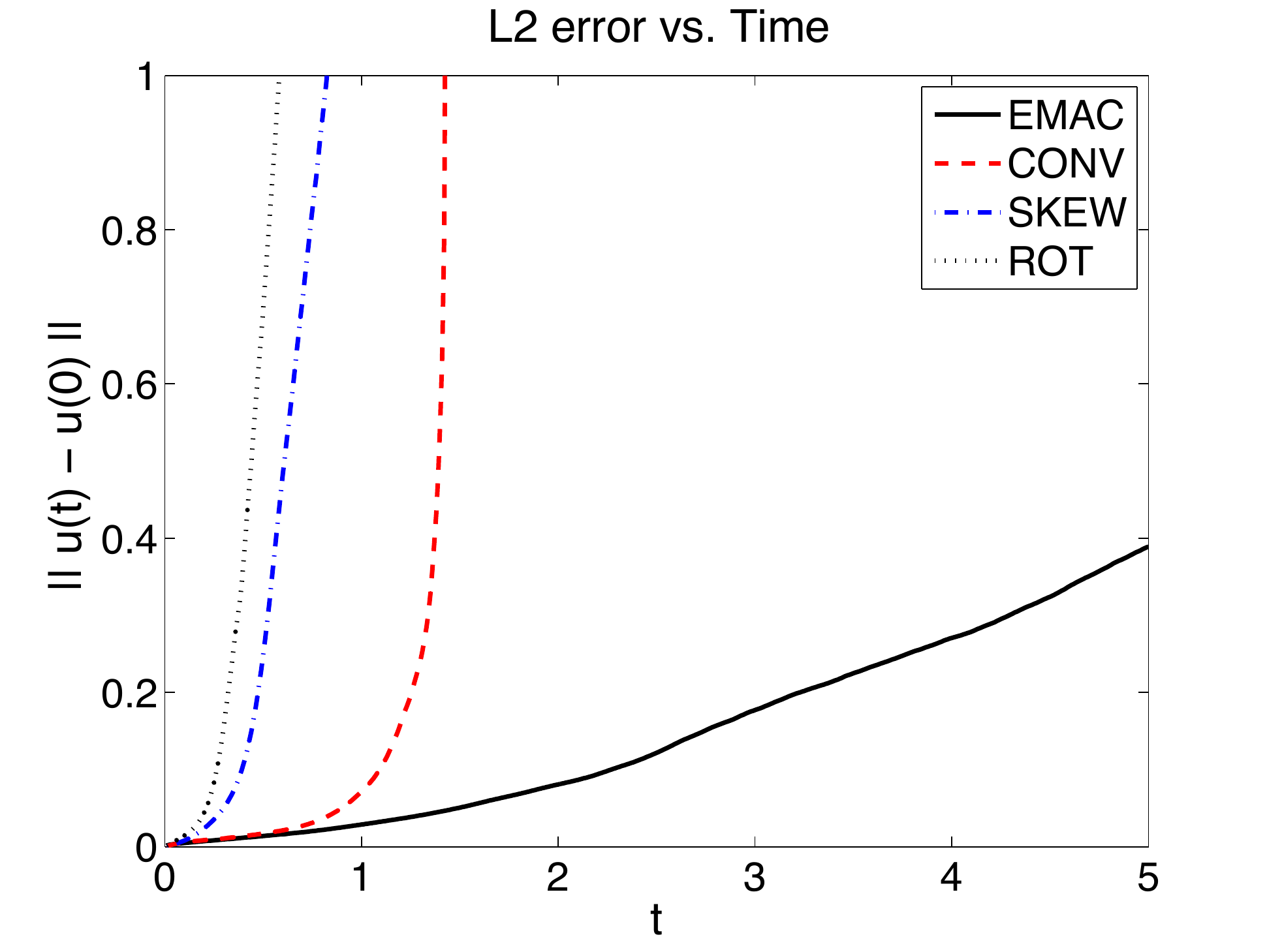}
\end{center}
\caption{\label{lattice1}
Shown above are plots of time versus energy, momentum, angular momentum, and $L^2(\Omega)$ velocity error, for EMAC, CONV, SKEW, and ROT formulations, for the lattice vortex problem.}
\end{figure}

We show results of energy, momentum, angular momentum, and $L^2(\Omega)$ error versus time, for each of the different formulations (using nonlinear solvers at each time step) in figure \ref{lattice1}.  Except for EMAC, the schemes all become unstable and blow up (exponential growth to $10^{16}$, at which point the simulation ends).

\begin{figure}[!h]
\begin{center}
\includegraphics[width=.49\textwidth,height=0.28\textwidth, viewport=0 0 550 420, clip]{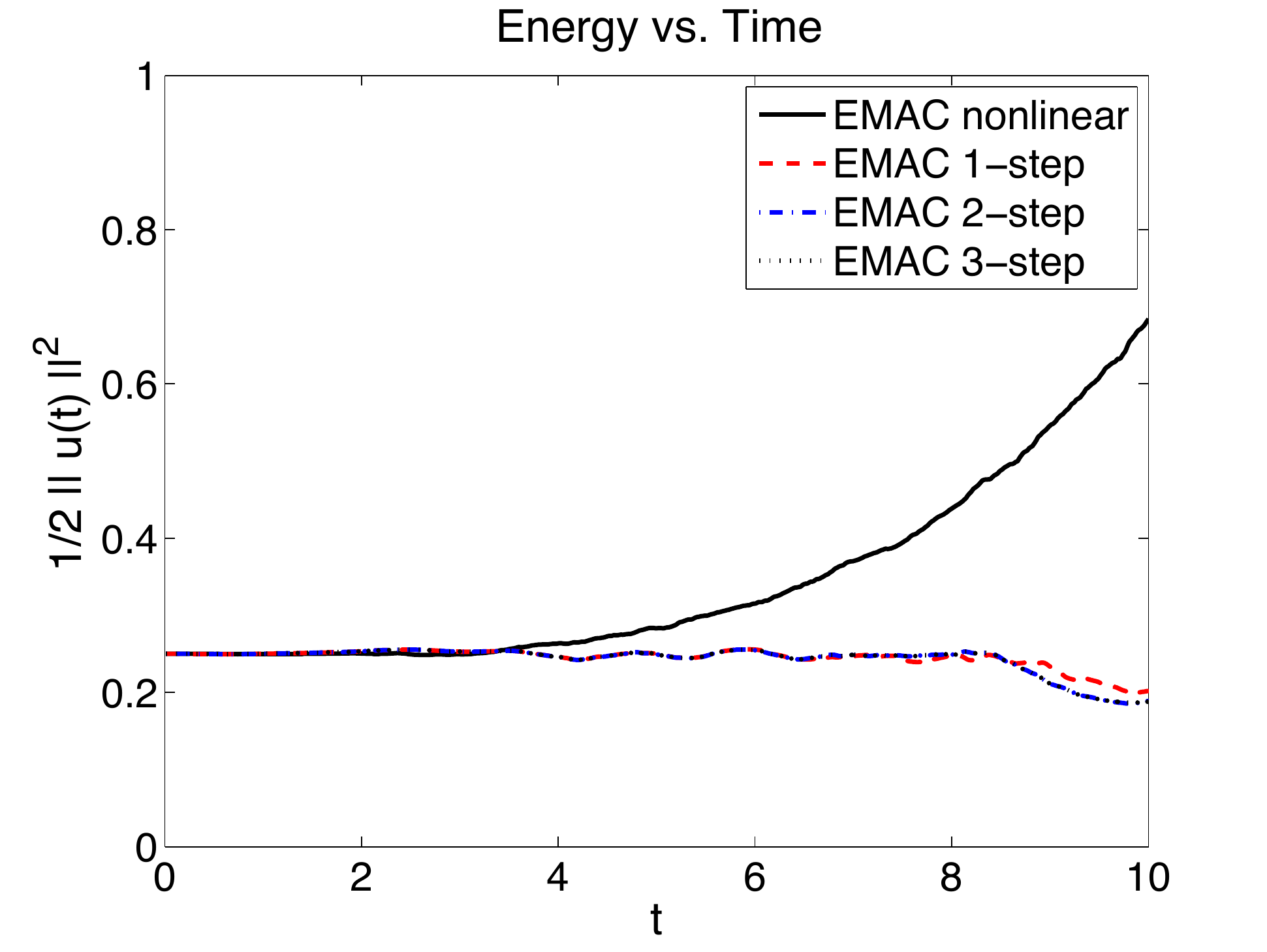}
\includegraphics[width=.49\textwidth,height=0.28\textwidth, viewport=0 0 550 420, clip]{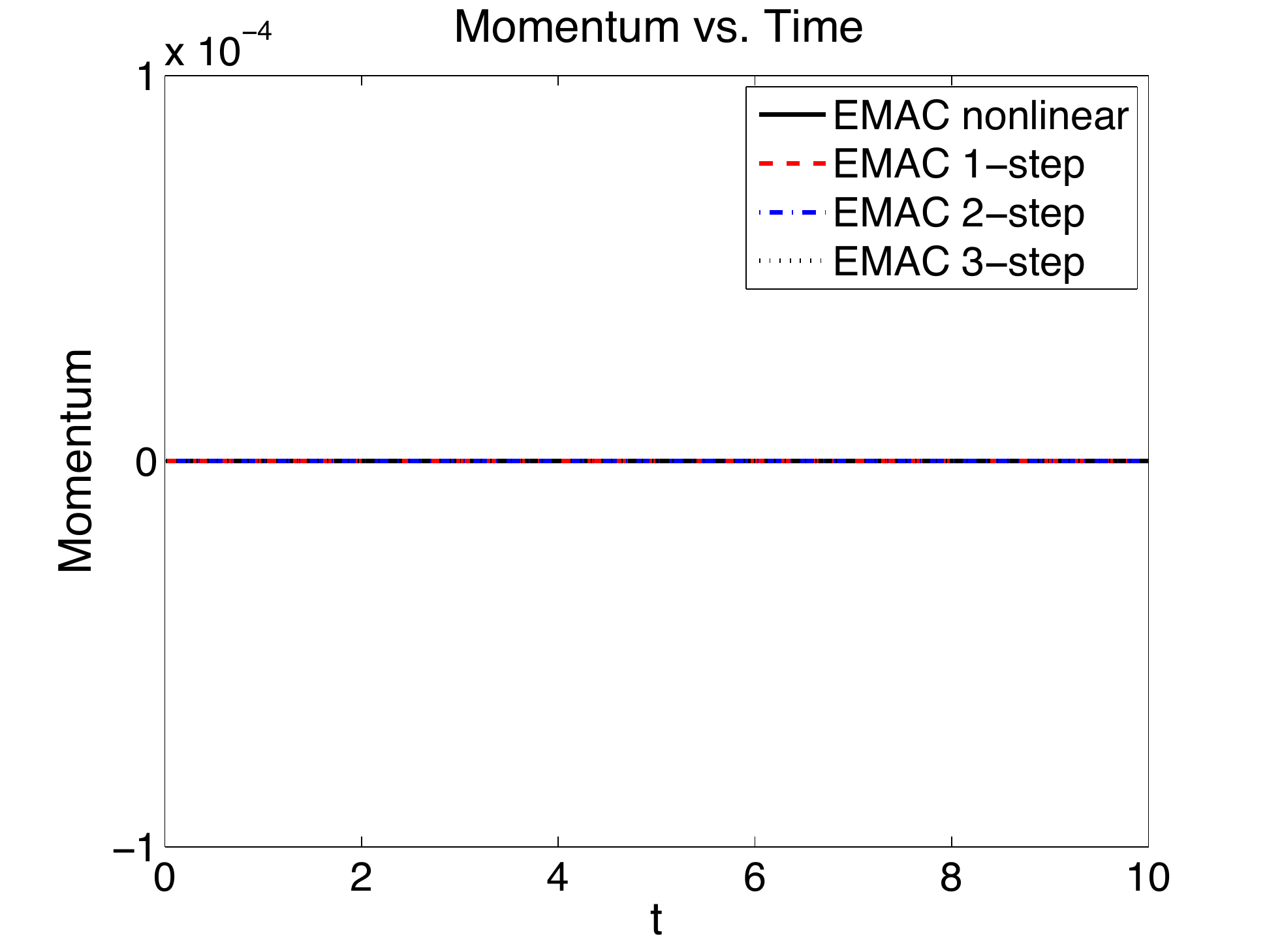}
\includegraphics[width=.49\textwidth,height=0.28\textwidth, viewport=0 0 550 420, clip]{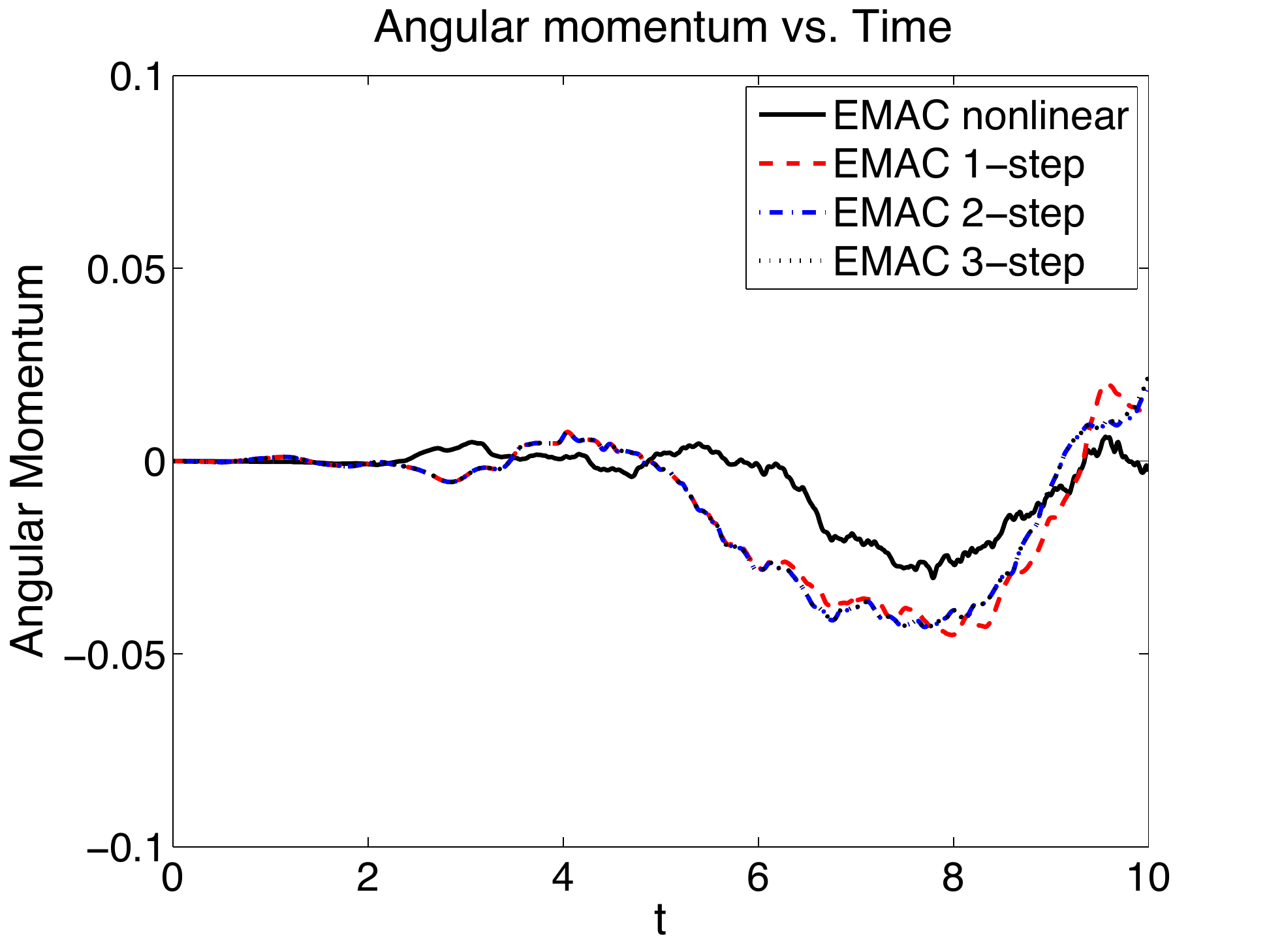}
\includegraphics[width=.49\textwidth,height=0.28\textwidth, viewport=0 0 550 420, clip]{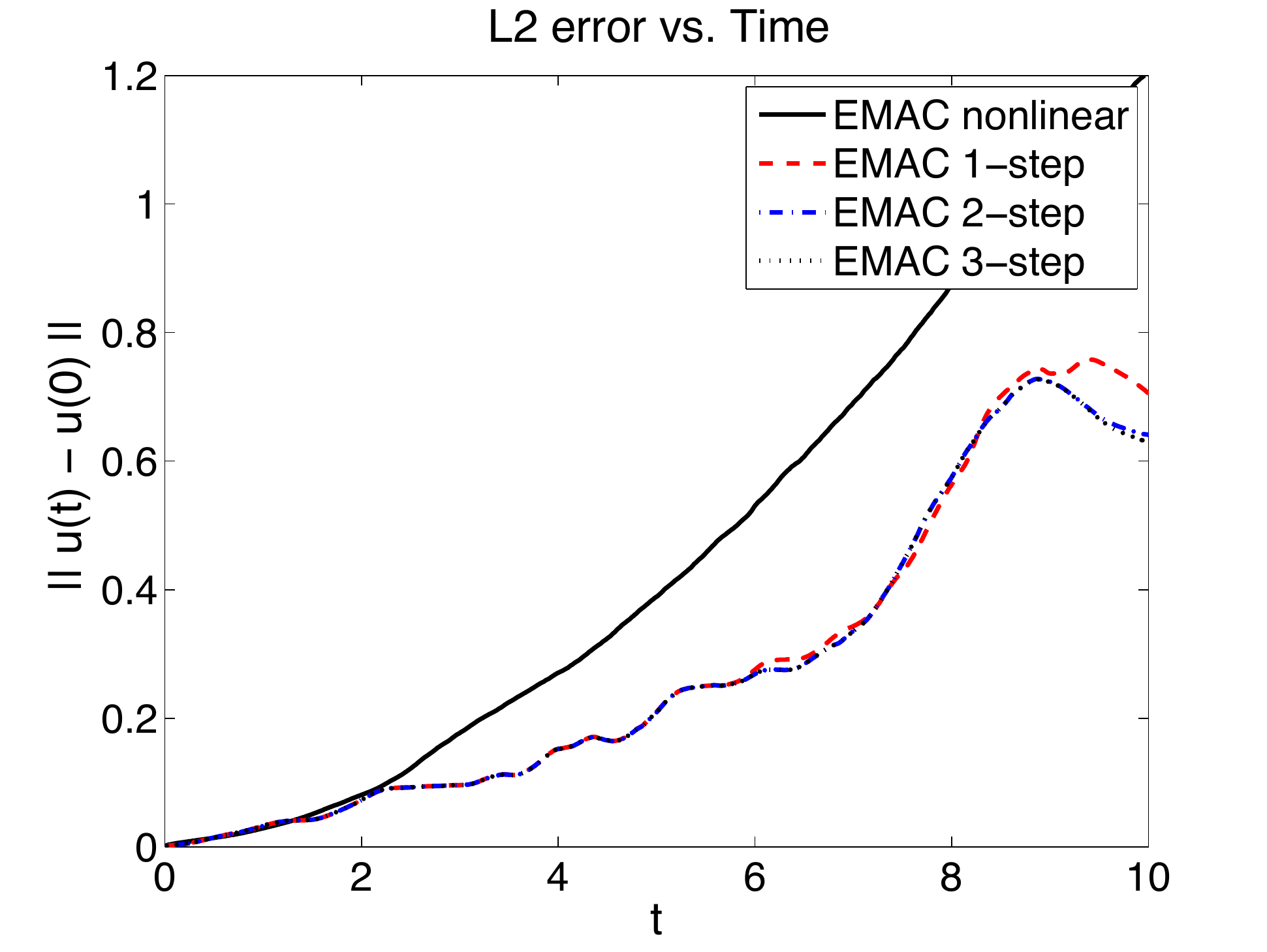}
\end{center}
\caption{\label{lattice2}
Shown above are plots of time versus energy, momentum, angular momentum, $L^2(\Omega)$ velocity error, for EMAC, CONV, SKEW, and ROT formulations, for the lattice vortex problem.}
\end{figure}

We show in figure \ref{lattice2} results up to t=10 for the same problem, but using EMAC with full nonlinear solve at each time step, and 1, 2, and 3 Newton steps at each time step.  We observe overall similar behavior for the linearizations of EMAC, compared to the fully nonlinear scheme; in fact, the linearizations are slightly better for error and energy, although there is no theoretical evidence to suggest this should be the case in general.

\subsection{3D flow past a circular cylinder}

This problem was first studied numerically for Reynolds number 20 and 100 by Schafer and Turek \cite{ref:schafer1996benchmark} in 1996. While the $Re=20$ problem gives a steady solution, higher Reynolds numbers lead to a time dependent solution. We consider the the case ``3D-3Z'' with the inlet being forced periodically in time leading
to a $0 \leq Re(t) \leq 100$.

\begin{figure}[!h]
\begin{center}
\includegraphics[width=.6\textwidth, clip]{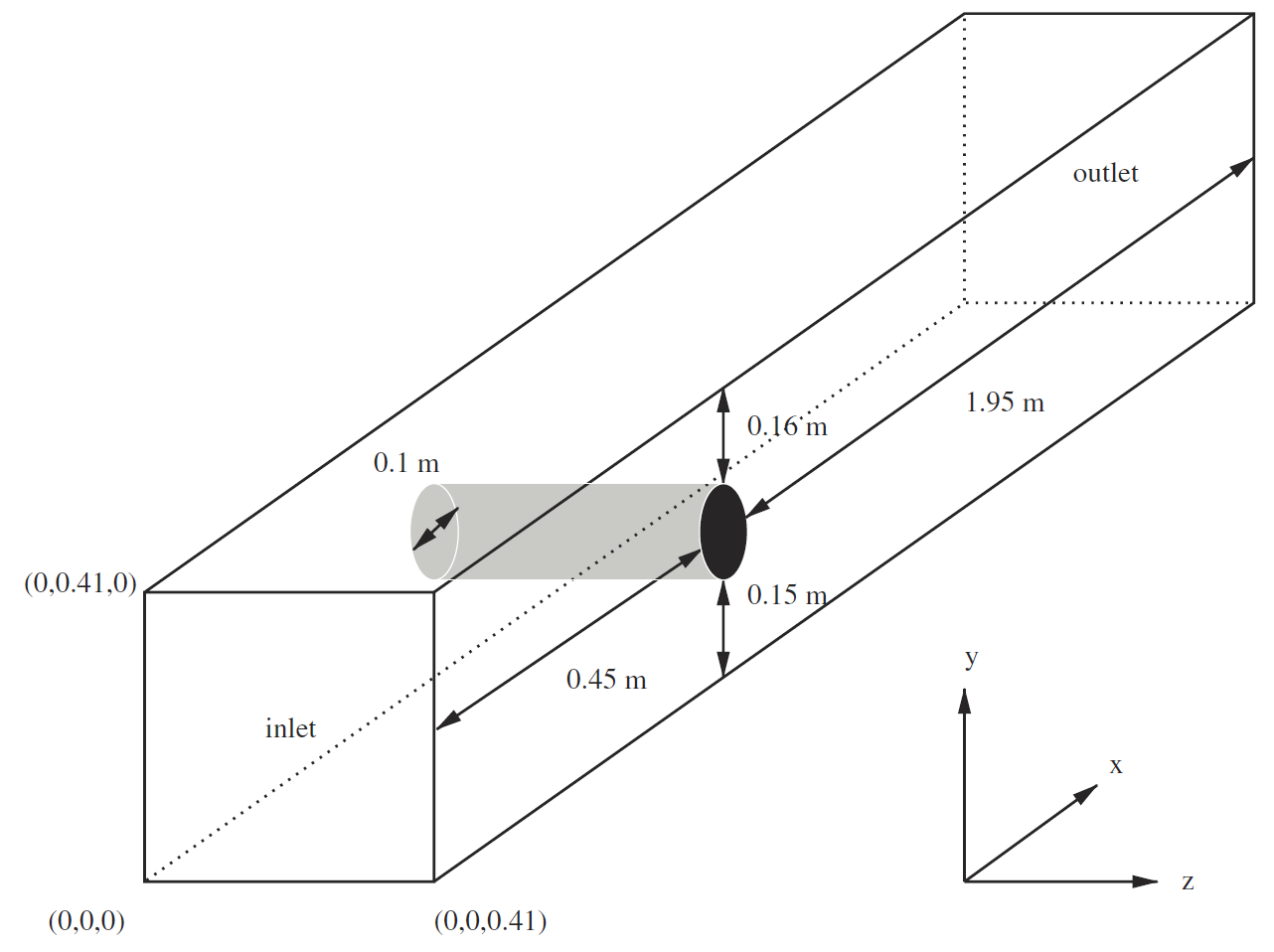}
\end{center}
\caption{\label{fig:Proj3:ChannelDiagram}
The channel with the circular cylinder.}
\end{figure}

The domain $\Omega$ is a 3D box with dimensions $0.41\times 0.41\times 2.5$m, and the obstacle is a circular cylinder with diameter $D=0.1$m. A diagram is shown in figure \ref{fig:Proj3:ChannelDiagram}.
Here, we denote $\Gamma_{walls}$ to be the bottom, left, right and top walls
of the channel and the boundary of the cylinder, $\Gamma_{in}$ to be the left boundary of the channel (inlet), and $\Gamma_{out}$ to be the right boundary (outlet).

The inlet flow profile is given for $0\le t \le 8$ (the times of interest) as
\[
		{\bf u}_x(0,y,z,t)|_{\Gamma_{in}} = \frac{16U_m yz(H-y)(H-z)}{H^4} \sin\left(\frac{\pi t}{8} \right) ,
\]
where $U_m = 2.25$ m/s, $H=0.41$m,  and we enforce no slip boundary conditions on the walls and cylinder, and use a zero traction condition at the outflow.

The computations are done in a parallel code using the deal.II library \cite{dealII84}. We use Trilinos \cite{trilinos}
for parallel linear algebra including ML as an algebraic multigrid preconditioner for the velocity block. For details
about the parallelization see \cite{BangerthBursteddeHeisterKronbichler11}, the block preconditioner based on
grad-div stabilization is the same as in \cite{CHOR17}. Also see \cite{heister_graddiv,zhao_liang_2017_484156} for more details.

We compute with $((Q_2)^3,Q_1)$ elements on a quadrilateral mesh refined heavily around the cylinder that provides 7.42 million total degrees of freedom.  The initial condition was taken to be the flow at rest, and a time step of $\Delta t=0.005$ was used to simulate the flow up to t=8.  A grad-div parameter of $\gamma=0.1$ was used for all simulations.
The runs take between 36 and 57 hours, using between 120 and 200 cores on the Palmetto cluster at Clemson University.

We report statistics for the different formulations in Table~\ref{table:3dcyl1} and the influence of the number
of nonlinear iterations in Table~\ref{table:3dcyl2}. The statistics are very similar between the different formulations
and it appears to be enough to do 2 Newton steps for the EMAC scheme, while a single Newton steps produces small variations in the statistics.  Interestingly, even with 7.42 million degrees of freedom, the simulations appear to not be fully resolved spatially, as max drag and lift from \cite{LS16,BMT12} predict max drag and lift to be $\sim$ 3.298 and 0.0028, respectively, using upwards of 90 million total degrees of freedom.  However, the test here was more about how much affect the linearization has on the answer, and Table~\ref{table:3dcyl2} shows it has very little effect.

\begin{table}
\centering
	\begin{tabular}{lrrrrrrr}
	Scheme & \multicolumn{1}{l}{max drag} & \multicolumn{1}{l}{min drag} & \multicolumn{1}{l}{min lift} & \multicolumn{1}{l}{max lift} & \multicolumn{1}{l}{$\Delta P_{max}$} & \multicolumn{1}{l}{$\Delta P_{min}$}\\  \hline
	CONS & 3.24648 & -0.170774 & -0.01022 & 0.002753 & 3.345390 & -0.103599 \\
	CONV & 3.24388 & -0.170789 & -0.01020 & 0.002753 & 3.345642 & -0.103568 \\
	EMAC & 3.25593 & -0.170797 & -0.01035 & 0.002753 & 3.351321 & -0.103541 \\
	ROT    & 3.24108 & -0.170805 & -0.01012 & 0.002753 & 3.336217 & -0.103574 \\
	SKEW & 3.24519 & -0.170781 & -0.01021 & 0.002753 & 3.345482 & -0.103583 \\
	\end{tabular}%
 	\caption{Statistics for the 5 formulations, where the nonlinear problem is fully resolved at each time step.}
\label{table:3dcyl1} 	
\end{table}

\begin{table}
\centering
	\begin{tabular}{lrrrrrrrr}
	Scheme (EMAC) & \multicolumn{1}{l}{max drag} & \multicolumn{1}{l}{min drag} & \multicolumn{1}{l}{min lift} & \multicolumn{1}{l}{max lift} & \multicolumn{1}{l}{$\Delta P_{max}$} & \multicolumn{1}{l}{$\Delta P_{min}$} \\ \hline
	Full nonlinear    & 3.25594 &     -0.17080 & -0.01035 & 0.002753 & 3.351321 & -0.103541 \\
	2-step Newton     & 3.25594 & -0.17080 & -0.01035 & 0.002753 & 3.351321 & -0.103541 \\
	1-step Newton    & 3.25962 &  -0.17087 & -0.00992 & 0.002762 & 3.352163 & -0.103532 \\
	\end{tabular}%
	\caption{EMAC statistics for 1, 2, and `as many as necessary' Newton steps take at each time steps.}
	\label{table:3dcyl2}
\end{table}


\section{Conclusions}

We have proposed and studied linearization methods in discretization of the EMAC formulation
of the NSE.  We have found that a skew-symmetric linearization, while energy stable, does not preserve
momentum or angular momentum, and moreover it performed poorly in our numerical tests.  However, the Newton
linearization was found to preserve both momentum and angular momentum, and although it does not conserve
energy exactly, we show that it is conserved up to the convergence of the Newton iteration, and also up to higher order
terms than the discretization error (for smooth solutions). Further, the EMAC Newton linearization performed very well in numerical tests: in all tests, if two steps of the EMAC Newton linearization were performed at each time step, the results were nearly identical to those of the nonlinear EMAC scheme.

\section{Acknowledgements}

Clemson University is acknowledged for generous allotment of compute time on Palmetto
cluster.

\bibliographystyle{plain}
\bibliography{graddiv}

\end{document}